\pdfoutput=1
\documentclass[final]{siamltex}
\usepackage{latexsym}
\usepackage{amsmath}
\usepackage{graphicx}
\usepackage{color}
\usepackage{graphicx}
\usepackage{amssymb}
\usepackage{amsmath}
\newtheorem{remark}{Remark}[section]
\newcommand{\vepsi}{\varepsilon}
\newcommand{\Del}{\Delta}
\newcommand{\Ome}{\Omega}
\newcommand{\p}{\partial}
\newcommand{\Div}{\mbox{div\,}}
\newcommand{\nab}{\nabla}

\title{Vanishing moment method and  moment solutions for second order fully 
nonlinear partial differential equations\footnote{This work is partially supported by the NSF grants DMS-0410266 and DMS-0710831.}}
\author{Xiaobing Feng\thanks{Department of Mathematics, 
The University of Tennessee, Knoxville, TN 37996, U.S.A. ({\tt xfeng@math.utk.edu}).}
\and Michael Neilan \thanks{Department of Mathematics, 
The University of Tennessee, Knoxville, TN 37996, U.S.A. ({\tt neilan@math.utk.edu}). }}
                                                                                
\begin{document}
\maketitle

\setcounter{page}{1}

\begin{abstract}
This paper concerns with numerical approximations of
solutions of second order fully nonlinear partial differential equations
(PDEs). A new notion of weak solutions, called moment solutions, is
introduced for second order fully nonlinear PDEs. Unlike
viscosity solutions, moment solutions are defined by a constructive
method, called vanishing moment method, hence, they can be
readily computed by existing numerical methods such as
finite difference, finite element, spectral Galerkin,
and discontinuous Galerkin methods with ``guaranteed" convergence.
The main idea of the proposed vanishing moment method is
to approximate a second order fully nonlinear PDE by
a higher order, in particular, a fourth order quasilinear PDE.
We show by various numerical experiments the viability
of the proposed vanishing moment method. All our numerical
experiments show the convergence of the
vanishing moment method, and they also show that
moment solutions coincide with viscosity solutions
whenever the latter exist.
\end{abstract}

\begin{keywords}
Fully nonlinear PDEs, Monge-Amp\`ere type equations,
 moment solutions, vanishing moment method, viscosity 
solutions, finite element method, mixed finite element method,
spectral and discontinuous Galerkin methods.
\end{keywords}

\begin{AMS}
65N30, 
65M60,  
35J60, 
35K55, 
53C45  
\end{AMS}

\section{Introduction}\label{sec-1}
Fully nonlinear PDEs are those PDEs which depend nonlinearly
on the highest order derivatives of unknown functions. 
Fully nonlinear PDEs arise from many areas in science and engineering
such as kinetic theory, materials science, differential
geometry, general relativity, optimal control, mass transportation,
image processing and computer vision, meteorology,
semigeostrophic fluid dynamics. They constitute the most difficult 
class of differential equations to analyze analytically and 
to approximate numerically, see \cite{Caffarelli_Cabre95,Gilbarg_Trudinger01,Fleming_Soner06,Caselles_Morel_Sbert98,Guan_Spruck04,McCann_Oberman04} and
references therein.

The general first order fully nonlinear PDE has the form 
\begin{equation}\label{eq1}
\quad F(\nabla u(x), u(x), x)=0 \qquad x\in \Ome\subset \mathbf{R}^n.
\end{equation}
The best known examples include Eikonal equation 
\[
|\nabla u(x)|=f(x)  \qquad x\in \Ome,
\]
and the general Hamilton-Jacobi equation 
\cite{Gilbarg_Trudinger01,Crandall_Ishii_Lions92}
\begin{equation*}
H(\nabla u(x))=0 \qquad x\in \Ome.
\end{equation*}
The general second order fully nonlinear PDE, which will be the focus of
this paper, takes the form
\begin{equation}\label{eq2}
\quad F(D^2 u(x),\nabla u(x), u(x), x)=0 \qquad x\in \Ome.
\end{equation}
where and throughout this paper $D^2u(x)$ denotes the
Hessian matrix of $u$ at $x$. The best known examples are the
Monge-Amp\`ere type equations \cite{Gilbarg_Trudinger01,Gutierrez01,Lieberman96}
\[
\mbox{det}(D^2u(x))=f(\nabla u(x),u(x),x) \qquad x\in \Ome,
\]
and the Bellman equations \cite{Gilbarg_Trudinger01,Fleming_Soner06}
\begin{equation}\label{bellman}
\displaystyle{\sup_{\theta\in \Theta} L_\theta(D^2 u,\nabla u, u, x)=0},
\end{equation}
where $\mbox{det}(D^2u(x))$ stands for the determinant of $D^2u(x)$ and
$L_\theta$ is a given family of second order linear differential
operators.

For the first order fully nonlinear PDEs, tremendous progresses
have been made in the past three decades. A revolutionary viscosity solution 
theory has been established (cf. \cite{Crandall_Lions83,Crandall_Evans_Lions84,Crandall_Ishii_Lions92,Fleming_Soner06}) and wealthy
amount of efficient and robust numerical methods and algorithms
have been developed and implemented (cf. \cite{Bryson_Levy03,Cockburn03,Crandall_Lions96,Lin_Tadmor00,Osher_Shu91,Zhang_Shu02,Zhao05}). However, for second order fully nonlinear PDEs,
the situation is strikingly different. On one hand, there have been
enormous advances in theoretical analysis 
in the past two decades after the introduction of the notion of
viscosity solutions by M. Crandall and  P. L. Lions in 1983 (cf.
\cite{Caffarelli03,Caffarelli_Cabre95,Caffarelli_Nirenberg_Spruck84,Crandall_Ishii_Lions92,Gilbarg_Trudinger01,Gutierrez01}).
On the other hand, in contrast to the success of the PDE analysis,
numerical solutions for general second order fully nonlinear PDEs
(except in the case of Bellman type PDEs, see below for details)
is mostly an untouched area, and computing viscosity solutions
of second order fully nonlinear PDEs has been
impracticable. There are several reasons for this lack of
progress. Firstly, the strong nonlinearity is an obvious one. Secondly,
the conditional uniqueness (i.e., uniqueness holds only in certain class
of functions) of solutions is difficult to handle numerically.
Lastly and most importantly, the notion of viscosity solutions, 
which is not variational, has no equivalence at the discrete level. 

To see the above points, let us consider the following model Dirichlet problem
for the Monge-Amp\`ere equation:
\begin{align}\label{ma}
\mbox{det}(D^2 u)&=f \qquad \mbox{in } \Ome, \\
u&=g \qquad \mbox{on }\p\Ome. \label{mab}
\end{align}
It is well-known that for non-strictly convex domain $\Ome$ the above problem
does not have classical solutions in general even $f$, $g$ and $\p\Ome$
are smooth (see \cite{Gilbarg_Trudinger01}). Classical
result of A. D. Aleksandrov states that the
Dirichlet problem with $f>0$ has a unique generalized solution in
the class of convex functions (cf. \cite{Aleksandrov61,Cheng_Yau77}). 
Major progress on analysis of problem \eqref{ma}-\eqref{mab}
has been made later by using the {\em viscosity solution} 
concept and machinery 
(cf. \cite{Caffarelli_Cabre95,Crandall_Ishii_Lions92,Gutierrez01}).
We recall that a {\em convex} function $u\in C^0(\overline{\Ome})$ 
satisfying $u=g$ on $\p \Ome$ is called a {\em viscosity subsolution 
(resp. viscosity supersolution)} of \eqref{ma} if for any $\varphi\in C^2$ 
there holds $\det(D^2\varphi(x_0))\leq f(x_0)$ (resp. 
$\det(D^2\varphi(x_0))\geq f(x_0)$) provided
that $u-\varphi$ has a local maximum (resp. a local minimum)
at $x_0\in \Ome$. $u\in C^0(\overline{\Ome})$ is called a {\em viscosity 
solution} if it is both a {\em viscosity subsolution} and 
a {\em viscosity supersolution}.

{\em First}, the reason to restrict the admissible set to be
the set of convex functions is that the Monge-Amp\`ere equation
is {\em elliptic} only in that set \cite{Gutierrez01}. 
It should be noted that in general the Dirichlet problem 
\eqref{ma}--\eqref{mab} may have other (nonconvex) solutions
besides the unique convex solution, multiple solutions are often expected
for the Monge-Amp\`ere type PDEs and for second order fully nonlinear
PDEs. It is easy to see that if one discretizes \eqref{ma}
straightforwardly using the finite difference method, one immediately
loses control on which solution the numerical scheme approximates even
assuming that the nonlinear discrete problem has solutions.
{\em Second}, the situation is even worse if one tries to formulate
a Galerkin type method (such as the finite element method and the spectral
Galerkin method), because there is no variational or weak formulation
to start with. In fact, this is clear from the definition
of {\em viscosity solutions}. It is not defined by the traditional
integration by parts approach, instead, it is defined by
a ``{\em differentiation by parts}" (a terminology coined by
L. C. Evans \cite{Evans90,Evans98}) approach. Although the 
``{\em differentiation
by parts}" approach has worked remarkably well for establishing
the viscosity solution theory for second order fully nonlinear PDEs in
the past two decades, it is extremely difficult (if all possible)
to mimic it at the discrete level.
{\em Third}, regardless which method is used, one can easily
envisage that the anticipated algebraic problem from
the discretization of a fully nonlinear PDE such as
the Monge-Amp\`ere equation must be very difficult to solve
due to the nonuniqueness of solutions and very strong nonlinearity.

Nevertheless, a few recent numerical attempts and results have been known 
in the literature. In \cite{Oliker_Prussner88} Oliker 
and Prussner proposed a finite difference scheme for 
computing Aleksandrov measure induced by $D^2u$ (and obtained 
the solution $u$ of \eqref{ma} as a by-product) in $2$-d.  
The scheme is extremely geometric and difficult to use and 
to generalize to other second order fully nonlinear PDEs. 
In \cite{Barles_Souganidis91} Barles and Souganidis
showed that any monotone, stable and consistent finite difference
scheme converges to the correct solution provided that there 
exists a comparison principle for the limiting equation.
Their result provides a guideline for constructing 
convergent finite difference methods although it did
not address how to construct such a scheme. Very recently, Oberman \cite{Oberman07}
was able to construct some wide stencil finite difference schemes
which fullfil the criterions listed in \cite{Barles_Souganidis91}
for the Monge-Amp\`ere type equations. 
In \cite{Baginski_Whitaker96} Baginski and
Whitaker proposed a finite difference scheme for Gauss curvature equation
(see \S \ref{sec-3}) in $2$-d by mimicking the unique continuation
method (used to prove existence of the PDE) at the discrete level.
Finally, in a series of papers 
\cite{Dean_Glowinski03,Dean_Glowinski04,Dean_Glowinski05,Dean_Glowinski06b} 
Dean and Glowinski proposed an augmented Lagrange multiplier method
and a least squares method for problem \eqref{ma}-\eqref{mab}
and the Pucci's equation (cf. \cite{Caffarelli_Cabre95,Gilbarg_Trudinger01})
in $2$-d by treating the Monge-Amp\`ere equation and Pucci's
equation as a constraint and using a variational criterion to select
a particular solution. Numerical experiments results were reported in
\cite{Oliker_Prussner88,Oberman07,Baginski_Whitaker96,Dean_Glowinski03,Dean_Glowinski04,Dean_Glowinski05,Dean_Glowinski06b}, however, convergence analysis was not addressed 
except in \cite{Oberman07}.

In addition, we like to remark that there
is a considerable amount of literature available on
using finite difference methods to approximate viscosity
solutions of second order fully nonlinear Bellman type
PDE \eqref{bellman} arisen from stochastic optimal control.
See \cite{Barles_Souganidis91,Barles_Jakobsen05,Jakobsen03,Krylov05}.
Due to the special nonlinearity of the Bellman type
PDEs, the approach used and the methods proposed in those papers
unfortunately could not be extended to other types of second
order fully nonlinear PDEs since the construction of those methods
critically relies on the linearity of the operators $L_\theta$.

The first goal of this paper is to introduce a new weak solution 
concept and a method to construct such a solution for second 
order fully nonlinear PDEs, in particular, for the Monge-Amp\`ere
type equations. These new weak solutions are called 
{\em moment solutions} and the method to construct such a 
moment solution is called {\em the vanishing moment method}.
The crux of this new method is that we approximate a second
order fully nonlinear PDE by a sequence of higher order (in particular, 
fourth order) {\em quasilinear} PDEs.  The limit of the solution sequence
of the higher order PDEs, if exists, is defined as a moment solution
of the original second order fully nonlinear PDE. Hence, moment solutions
are constructive by nature.  The second goal of this paper
is to present a number of numerical methods for computing 
moment solutions of second order fully nonlinear 
PDEs, and to present extensive numerical experiment results 
to demonstrate the convergence and effectiveness of the 
proposed vanishing moment methodology. Indeed, one of 
advantages of the vanishing moment method is that it allows
one to use wealthy amount of existing numerical methods and
algorithms as well as computer codes for fourth order linear 
and quasilinear PDEs to solve second order fully nonlinear PDEs.
The third and the last goal of this paper is to show using
numerical studies that the notion of moment solutions
generalizes the notion of viscosity solutions in the sense that
the former coincides with the later whenever the later exists.
These numerical studies indeed motivate us to give a 
rigorous convergence analysis of the vanishing moment method 
for the Monge-Amp\`ere equation in two spatial dimensions
\cite{feng06c}.

The remainder of the paper is organized as follows. In \S \ref{sec-2}, we
introduce the abstract framework of moment solutions and the vanishing 
moment method for general second order fully nonlinear PDEs.
In \S \ref{sec-4}, we propose two classes of numerical 
discretization methods and briefly discuss solution algorithms.
In \S \ref{sec-3}, we apply the abstract framework to several 
classes of second order fully nonlinear PDEs which include 
the Monge-Amp\`ere type equations, Pucci's extremal equations,
the infinite Laplace equation, and second order parabolic 
fully nonlinear PDEs. In \S \ref{sec-5}, 
we present many $2$-d and $3$-d numerical experiment 
results to demonstrate the convergence and effectiveness
of the vanishing moment methodology, and provide numerical evidences
of the agreement of moment solutions and viscosity solutions
whenever the latter exists. The paper is concluded by
a summary and some conclusions in \S 6.

\section{Vanishing moment method and the notion of moment solutions} 
\label{sec-2}

\subsection{Preliminaries} \label{sec-2.1}
Standard space notation will be adopted throughout this paper, 
we refer to \cite{Gilbarg_Trudinger01,Lieberman96} for 
their exact definitions. $\Ome$ denotes a generic
bounded domain in $\mathbf{R}^n$. $(\cdot,\cdot)$ and
$\langle\cdot, \cdot\rangle$ are used to denote the $L^2$-inner
products on $\Ome$ and on $\p\Ome$, respectively. We assume $n\geq 2$,
except in \S \ref{sec-4} and \S \ref{sec-5}, where we restrict
$n=2,3$ when we develop numerical methods and perform numerical experiments. 

Since the notion of viscosity solutions has been and will continue
to be referred many times, and it is closely related
to the notion of moment solutions to be described later in 
this paper, for readers' convenience, we briefly 
recall its definition and history here and refer to 
\cite{Caffarelli_Cabre95,Crandall_Ishii_Lions92,Crandall_Lions83,Evans98} 
for detailed discussions.

\begin{definition}\label{def1}
Suppose $F: \mathbf{R}^n\times \mathbf{R}\times \mathbf{R}^n 
\rightarrow \mathbf{R}$ is continuous (nonlinear)
function. 
\begin{enumerate}
\item[(i)] A function $u\in C^0(\Ome)$ is called a {\em viscosity subsolution}
of \eqref{eq1} if, for every $C^1$ function $\varphi=\varphi(x)$ such that
$u-\varphi$ has a local maximum at $x^0\in \Ome$, there holds
\[
F(\nabla \varphi(x^0),\varphi(x^0),x^0) \leq 0.
\]
\item[(ii)] A function $u\in C^0(\Ome)$ is called a {\em viscosity 
supersolution} of \eqref{eq1} if, for every $C^1$ function 
$\varphi=\varphi(x)$ such that $u-\varphi$ has a local minimum 
at $x^0\in \Ome$, there holds
\[
F(\nabla \varphi(x^0),\varphi(x^0),x^0) \geq 0.
\]
\item[(iii)] A function $u\in C^0(\Ome)$ is called a {\em viscosity 
solution} of \eqref{eq1} if it is both a viscosity subsolution and 
a viscosity supersolution.
\end{enumerate}
\end{definition}

It should be pointed out that the above definition is a modern
definition of viscosity solutions for \eqref{eq1}. It can be regarded as a  
``differentiation by parts" definition (cf. \cite{Evans98}). However,
viscosity solutions were first introduced differently by a vanishing viscosity 
procedure (cf. \cite{Crandall_Lions83}), that is, equation \eqref{eq1} 
is approximated by the second order {\em quasilinear} PDEs
\[
-\epsilon \Delta u^\epsilon + F(\nabla u^\epsilon,u^\epsilon,x) =0,
\]
and $\underset{\vepsi\rightarrow 0^+}{\lim} u^\vepsi$, if exists, is called
a viscosity solution of \eqref{eq1}. It was later proved that the 
two definitions are equivalent for equation \eqref{eq1} 
(cf. \cite{Crandall_Evans_Lions84}). 

Another important reason to favor the modern ``differentiation by parts" 
definition is that the definition and the notion of viscosity solutions 
can be readily extended to second order fully nonlinear PDEs. 

\begin{definition}\label{def2}
Suppose $F:\mathbf{R}^{n\times n}\times \mathbf{R}^n\times 
\mathbf{R}\times \mathbf{R}^n \rightarrow \mathbf{R}$ is 
continuous (nonlinear) function.
\begin{enumerate}
\item[(i)] A function $u\in C^0(\Ome)$ is called a {\em viscosity subsolution}
of \eqref{eq2} if, for every $C^2$ function $\varphi=\varphi(x)$ such that
$u-\varphi$ has a local maximum at $x^0\in \Ome$, there holds
\[
F(D^2\varphi(x^0),\nabla \varphi(x^0), \varphi(x^0), x^0) \leq 0.
\]
\item[(ii)] A function $u\in C^0(\Ome)$ is called a {\em viscosity
supersolution} of \eqref{eq2} if, for every $C^2$ function
$\varphi=\varphi(x)$ such that $u-\varphi$ has a local minimum
at $x^0\in \Ome$, there holds
\[
F(D^2\varphi(x^0),\nabla \varphi(x^0), \varphi(x^0), x^0) \geq 0.
\]
\item[(iii)] A function $u\in C^0(\Ome)$ is called a {\em viscosity
solution} of \eqref{eq2} if it is both a viscosity subsolution and
a viscosity supersolution.
\end{enumerate}
\end{definition}

As it is known now, a successful theory of viscosity solutions has been 
established for second order fully nonlinear PDEs in the past two decades 
(cf. \cite{Caffarelli_Cabre95,Crandall_Ishii_Lions92,Gutierrez01}). 
On the other hand, it should be noted that the phrase ``viscosity solution" 
loses its original meaning in this theory since it has nothing to do 
with the vanishing viscosity method in the case of second order fully 
nonlinear PDEs. We recall that to establish the existence of 
viscosity solutions the technique used to substitute for the vanishing 
viscosity method in the theory is the classical Perron's method
(cf. \cite{Crandall_Ishii_Lions92,Caffarelli_Cabre95}). 
To the best of our knowledge, viscosity solutions of second order 
fully nonlinear PDEs were never defined and/or constructed by a 
limiting process like one described above for the Hamilton-Jacobi equation.

\subsection{General framework of the vanishing moment method}\label{sec-2.2}

For the reasons and difficulties explained in \S \ref{sec-1}, as
far as we can see, it is unlikely (at least very difficult if 
all possible) that one can {\em directly} 
approximate viscosity solutions of general second order fully nonlinear
PDEs such as Monge-Amp\`ere type equations using any available
numerical methodology (finite difference method, finite element
method, spectral method, meshless method etc.).  From computational
point of view, the notion of viscosity solutions is a ``bad" notion
for second order fully nonlinear PDEs because it is not constructive nor
variational, so one has no handle on how to compute such a solution.

In searching for a ``better" notion of weak solutions for second order
fully nonlinear PDEs, we are inspired by the following simple but 
crucial observation: {\em the essence of the vanishing viscosity
method for the Hamilton-Jacobi equation and the original notion
of viscosity solutions is to approximate a lower order fully nonlinear
PDE by a sequence of higher order quasilinear PDEs}.
This observation then suggests us to apply the above principle
to second order fully nonlinear PDE \eqref{eq2}, this is exactly 
what we are going to do in this paper. That is, we approximate 
equation \eqref{eq2} 
by the following higher order quasilinear PDEs:
\begin{equation}\label{eq3}
G_\vepsi(D^r u^\vepsi) + F(D^2u^\vepsi,\nabla u^\vepsi,x)=0
\qquad (r\geq 3,\, \vepsi>0),
\end{equation}
where $\{G_\vepsi\}$ is a family of suitably chosen linear or quasilinear 
differential operators of order $r$. The above approximation 
then naturally leads to the next definition.

\begin{definition}\label{def3}
Suppose that $u^\vepsi$ solves \eqref{eq3} for each $\vepsi>0$, we call 
$\underset{\vepsi\rightarrow 0^+}{\lim} u^\vepsi$ a {\em moment solution} of 
\eqref{eq2} provided that the limit exists. We also call this limiting 
process {\em the vanishing moment method}.
\end{definition}

Clearly, the above definition is a loose definition since 
the operator $G_\vepsi$ is not specified, nor is the 
meaning of the limit, but they will become clear later in this section.
We note that the reason to use the terminology ``moment solution"
will also be explained later in this section, and the notion of moment 
solutions and the vanishing moment method are clearly in
the spirit of the (original) notion of viscosity solution
and the vanishing viscosity method \cite{Crandall_Lions83}. 

To establish a complete theory of moment solutions and
vanishing moment method for second order fully nonlinear PDEs,
there are many issues we must address. For instance,
\begin{itemize}
\item How to choose the operator $G_\vepsi$ ?
\item What additional boundary condition(s) should $u^\vepsi$ satisfy?
\item Does the limit $\underset{\vepsi\rightarrow 0^+}{\lim} u^\vepsi$
always exist? If it does, what is the rate of convergence?
\item How do moment solutions relate to viscosity solutions?
\item How to solve \eqref{eq3} numerically?
\item Error estimates, nonlinear solvers, computer implementations.
\end{itemize}
As expected, we do not have answers for all the questions now, 
nor do we intend to address all of them in this paper. Instead,
the focus of this paper is to develop the framework for 
moment solutions and the vanishing viscosity method, and to present
numerical evidences to show effectiveness of the method and to justify 
the proposed approach. On the other hand, we do plan to address
all theoretical issues in forthcoming papers \cite{feng06c,Neilan06}.  

Regarding to the first issue, although the choices for $G_\vepsi$ 
are abundant and flexible, the following are some guidelines for 
choosing a good operator $G_\vepsi$.
\begin{enumerate} 
\item[(a)] $G_\vepsi$ must be a linear or quasilinear operator. 
\item[(b)] $G_\vepsi\rightarrow 0$ in some reasonable sense 
as $\vepsi\rightarrow 0^+$. 
\item[(c)] $G_\vepsi(D^r u)$ is better to be elliptic, in particular, when 
PDE \eqref{eq2} is elliptic. 
\item[(d)] Equation \eqref{eq3} should be relatively easy 
to solve numerically.
\end{enumerate}

Since an elliptic operator is necessarily of even order, so 
guideline (c) above implies that $r$ must be an even number 
in \eqref{eq3}. Hence, the lowest order of equation \eqref{eq3} 
is $r=4$. When talking about fourth order elliptic operators, 
the biharmonic operator stands out immediately. So we let
\begin{equation*}
G_\vepsi(D^4 v):=-\vepsi\Delta^2 v,
\end{equation*}
then equation \eqref{eq3} becomes
\begin{equation}\label{ePDE}
-\vepsi\Delta^2 u^\vepsi + F(D^2u^\vepsi,\nabla u^\vepsi,x)=0.
\end{equation}
After the differential operator $G_\vepsi$ is chosen, next we need to
take care the boundary conditions. Here we only consider Dirichlet
problem for \eqref{eq2}. Suppose that 
\begin{equation}\label{bc1}
u=g\qquad\mbox{on}\quad \p\Ome, 
\end{equation}
it is obvious that we need to impose 
\begin{equation}\label{bc2}
u^\vepsi = g \quad\mbox{or}\quad u^\vepsi \approx g 
\qquad\mbox{on}\quad \p\Ome.
\end{equation}
Moreover, since \eqref{ePDE} is a fourth order PDE,
in order to uniquely determine $u^\vepsi$ we need to impose an
additional boundary condition for $u^\vepsi$. Mathematically, 
many boundary conditions can be used for this purpose. Physically, 
any additional boundary condition will introduce a ``boundary 
layer", so a better choice would be one which minimizes the boundary
layer. Here we proposed to use one of following three boundary conditions
\begin{equation}\label{bc3}
\Del u^\vepsi =\vepsi^2 \qquad\mbox{or}\qquad
\frac{\p \Del u^\vepsi}{\p \mathbf{n}}=\vepsi^2
\qquad\mbox{or}\qquad D^2u^\vepsi\mathbf{n}\cdot \mathbf{n} =\vepsi^2
\qquad\mbox{on}\quad \p\Ome.
\end{equation}
In particular, the first two boundary conditions, which are natural
boundary conditions, have an advantage in PDE convergence analysis 
\cite{feng06c,Neilan06}. Another valid boundary condition is the 
Neumann boundary condition $\frac{\p u^\vepsi}{\p \mathbf{n}}=\vepsi^2$
on $\p\Ome$. But since this is an essential boundary condition,
it produces a larger boundary layer than the above three boundary 
conditions. The rationale for picking the above  boundary conditions
is that we implicitly impose an extra boundary condition 
$\vepsi^m \Del u^\vepsi + u^\vepsi =g+ \vepsi^{m+2}$ 
on $\p\Ome$ for equation \eqref{ePDE}, which is a higher order
perturbation of the original Dirichlet boundary condition $u^\vepsi=g$
on $\p\Ome$. Intuitively, we expect and hope that the extra boundary 
condition converges to the original Dirichlet boundary condition as 
$\vepsi$ tends to zero. We note that $m$ can be any positive integer, 
and power $2$ is used for convenience and it can be replaced by 
any positive integer.

We now remark that when $n=2$ in mechanical applications
$u^\vepsi$ often stands for the vertical displacement of a plate
and $D^2u^\vepsi$ is the {\em moment tensor}, and in the weak formulation,
the biharmonic term becomes $-\vepsi (D^2u^\vepsi, D^2 v)$ which
should vanish as $\vepsi\rightarrow 0^+$.  This is the very reason 
why we call $\underset{\vepsi\rightarrow 0^+}{\lim} u^\vepsi$, if exists, 
a moment solution and call the limiting process the vanishing moment method.

In summary, we propose to approximate the second order fully nonlinear
Dirichlet problem \eqref{eq2},\eqref{bc1} by the fourth order 
quasilinear boundary value problems \eqref{ePDE},\eqref{bc2},
\eqref{bc3}. Since we expect $u^\vepsi\in W^{m,p}(\Ome)$ 
for $m\geq 2, p\geq 2$,  so the convergence 
$\underset{\vepsi\rightarrow 0^+}{\lim} u^\vepsi$ in Definition
\ref{def3} can be understood in $H^2$-topology or in $H^1$-topology
or even in $L^2$-topology. To distinguish these different limits,
we introduce the following refined definition of Definition \ref{def3}.

\begin{definition}\label{def4}
Suppose that $u^\vepsi\in H^2(\Ome)$ solves problem \eqref{ePDE},
\eqref{bc2}, \eqref{bc3}. $\underset{\vepsi\rightarrow 0^+}{\lim} u^\vepsi$
is called respectively a {\em sub-weak, weak} and {\em strong moment solution} 
to problem \eqref{eq2},\eqref{bc1} if the convergence 
holds in $L^2$-, $H^1$- and $H^2$-topology. 
\end{definition}

\begin{remark}
Since sub-weak and weak moment solutions do not have second order
weak derivatives, they are very hard (if all possible) to identify.
On the other hand, since strong moment solutions do have  second order
weak derivatives, naturally they are expected to satisfy the PDE 
\eqref{eq2} almost everywhere in $\Ome$ and to fulfill 
the boundary condition \eqref{bc1} pointwise on $\p\Ome$
(cf. \cite{feng06c,Neilan06}).
\end{remark}

\section{Discretization and solution methods} \label{sec-4}

The vanishing moment method reduces the problem of solving 
\eqref{eq2},\eqref{bc1} to a problem of solving 
\eqref{ePDE},\eqref{bc2}$_1$,\eqref{bc3}$_1$ for each
fixed $\vepsi>0$. Since \eqref{ePDE} is a nonlinear biharmonic
equation, one can use any of wealthy amount of existing 
numerical methods for biharmonic problems to discretize the equation. 
Although other types of numerical methods are applicable, 
here we focus on Galerkin
type methods such as finite element methods, mixed finite
element methods, discontinuous and spectral Galerkin methods
\cite{Ciarlet02,Brezzi_Fortin91,Cockburn_Karniadakis_Shu00,Bernardi_Maday97}.
Throughout this section, we assume $n=2,3$.

\subsection{Finite element methods in $2$-d} \label{sec-4.1}

In the two-dimensional case many finite element methods, such as
confirming Argyris, Bell, Bogner--Fox--Schmit and Hsieh--Clough--Tocher 
elements and nonconforming Adini, Morley, and Zienkiewicz elements,
were extensively developed in 60's and 70's for the biharmonic problems. 
A beautiful theory of plate finite element methods was also established   
(cf. \cite{Ciarlet02}). Naturally, one would want to solve 
problem \eqref{ePDE},\eqref{bc2}$_1$,\eqref{bc3}$_1$ by using and adapting
these well-known plate finite element methods. That is exactly 
what we are going to do next. For the sake of presentation clarity,
here we only discuss the confirming finite element methods,
and refer to \cite{Neilan06} for a detailed development of nonconfirming
finite element methods for problem \eqref{ePDE},\eqref{bc2}$_1$,
\eqref{bc3}$_1$.

The variational formulation for \eqref{ePDE},\eqref{bc2}$_1$,
\eqref{bc3}$_1$ is defined as: Find $u^\vepsi\in H^2(\Ome)$
with $u^\vepsi=g$ a.e. on $\p\Ome$ such that for any 
$v\in  H^2(\Ome)\cap H^1_0(\Ome)$ there holds
\begin{align} \label{weakform}
-\vepsi\bigl( \Del u^\vepsi, \Del v\bigr) 
+ \bigl( F(D^2 u^\vepsi,\nabla u^\vepsi,u^\vepsi,x), v\bigr) 
&= -\Bigl\langle \vepsi^3, \frac{\p v}{\p \mathbf{n}} \Bigr\rangle. 
\end{align}
Let $\mathcal{T}_h$ be a quasiuniform triangular or rectangular mesh 
with mesh size $h\in (0,1)$ for the domain $\Ome\subset \mathbf{R}^2$. 
Let $U^h_g\subset H^2(\Ome)$ 
denote one of confirming finite element spaces (as mentioned above)
whose functions take the boundary value $g$ at all nodes 
on $\p\Ome$. Then our finite element method is defined as: 
Find $u^\vepsi_h\in U^h_g$ such that 
\begin{align} \label{fem}
-\vepsi\bigl( \Del u^\vepsi_h, \Del v_h\bigr) 
+ \bigl( F(D^2 u^\vepsi_h,\nabla u^\vepsi_h,u^\vepsi_h,x), v_h\bigr) 
&= -\Bigl\langle \vepsi^3, \frac{\p v_h}{\p \mathbf{n}} \Bigr\rangle,
\quad \forall v_h\in  U^h_0.
\end{align}

In \S \ref{sec-5}, we shall present several numerical experiment results 
for the Monge-Amp\'ere type equations to show the excellent performance 
of the Argyis finite element method. Convergence and error analysis of 
the above scheme and other finite element schemes will be presented 
in forthcoming papers (also see \cite{Neilan06}).

\subsection{Mixed finite element methods in $2$-d and $3$-d} \label{sec-4.2}

Along with the theory of plate finite element methods, another beautiful
theory of mixed finite element methods was also extensively 
developed in '70s and '80s for the biharmonic problems in $2$-d 
(cf. \cite{Brezzi_Fortin91,Ciarlet02,Falk_Osborn80}).  It is interesting 
to point out that all these $2$-d mixed finite element methods can be 
easily generalized to solving $3$-d biharmonic problems and 
general fourth order quasilinear PDEs 
(cf. \cite{Elloitt_French_Miller89,XA2,XA3}).

Because the Hessian matrix $D^2 u^\vepsi$ appears in \eqref{ePDE}
in a nonlinear fashion,  to design a mixed method we are ``forced" to
introduce $\sigma^\vepsi:=D^2 u^\vepsi$ (not $v^\vepsi:=\Del u^\vepsi$
alone) as additional variables so the mixed method simultaneously
seeks $u^\vepsi$ and $\sigma^\vepsi$. This observation then excludes
the usage of the popular family of Ciarlet-Raviart mixed finite 
element methods (originally designed for the biharmonic problems) 
\cite{Ciarlet02,Ciarlet_Raviart74},
on the other hand, the observation suggests to try Hermann-Miyoshi
mixed elements 
\cite{Falk_Osborn80,Hermann67,Miyoshi73,Miyoshi76,Oukit_Pierre96} 
and Hermann-Johnson mixed elements
\cite{Falk_Osborn80,Hermann67,Johnson73} both use $\sigma^\vepsi$ as
additional variables.  

To define Hermann-Miyoshi type mixed finite element methods,
we first derive the following mixed variational formulation
for problem \eqref{ePDE},\eqref{bc2}$_1$,\eqref{bc3}$_3$: 
Find $(u^\vepsi,\sigma^\vepsi)\in V_g\times W_\vepsi$ such that
\begin{alignat}{2}\label{mixed1}
&(\sigma^\vepsi,\mu)+(\nabla u^\vepsi, \Div\mu) =\sum_{i=1}^{n-1}
\Bigl\langle \frac{\p g}{\p \mathbf{\tau}_i},
\mu \mathbf{n}\cdot \mathbf{\tau}_i \Bigr\rangle &&\qquad
\forall \mu\in W_0,\\
&\vepsi(\Div \sigma^\vepsi, \nabla v) 
+ \bigl( F(\sigma^\vepsi,\nabla u^\vepsi, u^\vepsi,x), v \bigr) 
= (f, v) &&\qquad \forall v\in V_0, \label{mixed2}
\end{alignat}
where $\mathbf{\tau}_i, i=1,2,\cdots, (n-1)$ denote the 
$(n-1)$ tangential directions at each point on $\p\Ome$, 
$\frac{\p g}{\p \mathbf{\tau}_i}$ denotes the tangential
derivative of $g$ along $\mathbf{\tau}_i$, and
\begin{align*}
V_g:&=\bigl\{ v\in H^1(\Ome);\, v|_{\p\Ome}=g \bigr\},\quad
V_0:=\bigl\{ v\in H^1(\Ome);\, v|_{\p\Ome}=0 \bigr\},\\
W_\vepsi:&=\bigl\{\mu \in [H^1(\Ome)]^{n\times n};\,\mu_{ij}=\mu_{ji}, 
\mu\mathbf{n}\cdot \mathbf{n}|_{\p\Ome}=\vepsi^2 \bigr\},\\
W_0:&=\bigl\{\mu \in [H^1(\Ome)]^{n\times n};\,\mu_{ij}=\mu_{ji},
\mu\mathbf{n}\cdot \mathbf{n}|_{\p\Ome}=0 \bigr\}.
\end{align*}

Let $\mathcal{T}_h$ be a quasiuniform triangular or rectangular mesh
if $n=2$ and be a quasiuniform tetrahedronal or $3$-d rectangular mesh
if $n=3$ for the domain $\Ome$. Let
$V^h\subset H^1(\Ome)$ be the Lagrange finite element space consisting of 
continuous piecewise polynomials of degree $k (\geq 2)$ associated with the
mesh $\mathcal{T}_h$. Let
\[
V^h_g:=V^h\cap V_g,  V^h_0:=V^h\cap V_0,\,\,
W^h_\vepsi:= [V^h]^{n\times n} \cap W_\vepsi,
W^h_0:= [V^h]^{n\times n} \cap W_0.
\]
Based on the variational formulation \eqref{mixed1}--\eqref{mixed2}
we define our (Hermann-Miyoshi type) mixed finite element methods 
as follows: Find $(u^\vepsi_h,\sigma^\vepsi_h)\in V^h_g\times W^h_\vepsi$ 
such that
\begin{alignat}{2}\label{mixed1a}
&(\sigma^\vepsi_h,\mu_h)+(\nabla u^\vepsi_h, \Div\mu_h)=\sum_{i=1}^{n-1}
\Bigl\langle \frac{\p g}{\p \mathbf{\tau}_i},
\mu_h \mathbf{n}\cdot \mathbf{\tau}_i \Bigr\rangle &&\quad
\forall \mu_h\in W^h_0,\\
&\vepsi(\Div \sigma^\vepsi_h, \nabla v_h) 
+\bigl( F(\sigma^\vepsi_h,\nabla u^\vepsi_h, u^\vepsi_h,x), v_h \bigr)
= (f, v_h) &&\quad \forall v_h\in V^h_0. \label{mixed2a}
\end{alignat}

Similarly, we can define variants of the above scheme as those proposed
in \cite{Oukit_Pierre96} as well as Hermann-Johnson type mixed methods. 
In \S \ref{sec-5}, we shall present several numerical experiment results
for the above scheme applying to the Monge-Amp\'ere type equations. 
Convergence and error analysis of the above scheme and 
other mixed finite element schemes will be presented in
a forthcoming paper (also see \cite{Neilan06}).

\begin{remark}
Besides the finite element and mixed finite element discretization methods,
one can also approximate problem \eqref{ePDE},\eqref{bc2}$_1$,\eqref{bc3}$_3$
by discontinuous Galerkin methods \cite{Arnold_Brezzi_Cockburn_Marini02,Cockburn_Kanschat_Schotzau05,Cockburn_Karniadakis_Shu00,Girault_Riviere_Wheeler05,fk05,fk06,Mozolevski_Suli03} and spectral Galerkin methods \cite{Bernardi_Maday97,Canuto_Hussaini_Quarteroni_Zang88,Shen94}. It should be pointed out that
these methods are dimenson-independent, hence, can be used in both 
$2$-d and $3$-d cases. We refer to \cite{Neilan06} for a detailed exposition. 
\end{remark}

\subsection{Remarks on second order fully nonlinear parabolic equations} 
\label{sec-4.5}

By adopting the method of line approach, generalizations of the 
numerical methods discussed in previous subsections to the corresponding 
parabolic equations \eqref{parabolic2} and \eqref{pma2} are standard 
(cf. \cite{Elloitt_French_Miller89,fk06} and references therein). 
Assuming that an implicit 
time stepping method such as the backward Euler and the 
Crank-Nicolson schemes will be used for time discretization,
then at each time step we only need to solve a fully nonlinear 
elliptic equation of the form \eqref{ePDE}. As
a result, all numerical methods discussed in 
\S \ref{sec-4.1}--\S \ref{sec-4.2} immediately apply.  
On the other hand, it should be pointed out that
the convergence and error analysis of all fully discrete
schemes are expected to be harder, in particular,
establishing error estimates which depend on $\vepsi^{-1}$
{\em polynomially} instead of {\em exponentially} will be very
challenging (cf. \cite{XA1,XA2,XA3,XA4,fw05,fw06}).

\subsection{Remarks on nonlinear solvers and preconditioning} \label{sec-4.6}

After equations \eqref{ePDE} and \eqref{parabolic2} are
discretized by any of above discretization methods, we get a 
strong nonlinear algebraic system
to solve. To the end, one has to use one or another
iterative solution method to do the job. 
In all numerical experiments to be given in \S\ref{sec-5},
we use preconditioned Newton iterative methods
as our nonlinear solvers. A few fixed point iterations 
might be needed to generate initial guess for Newton type 
iterative methods. Another strategy which we are currently
investigating is the following ``multi-resolution" strategy: 
first compute a numerical solution using a relatively 
large $\vepsi$, then use the computed solution as an 
initial guess for the Newton method at the finer resolution $\vepsi$.
Regarding to preconditioning, we use the simple
ILU preconditioner in all simulations of \S\ref{sec-5}. 
We plan to use more sophisticate multigrid and Schwarz (or domain
decomposition) preconditioners when the amount of computations
becomes intensive and large in $3$-d. 

\section{Applications} \label{sec-3}

In this section, we shall apply the vanishing moment methodology
outlined in the previous section to several classes of specific second
order fully nonlinear PDEs.

\subsection{Monge-Amp\`ere type equations} \label{sec-3.1}
Monge-Amp\`ere type equations refer to a class of second order
fully nonlinear PDEs of the form (cf. \cite{Caffarelli_Cabre95,Caffarelli_Nirenberg_Spruck84,Cheng_Yau77,Gilbarg_Trudinger01,Gutierrez01})
\begin{equation} \label{ma1}
F(D^2u^0, Du^0, u^0,x):=\mbox{det}(D^2 u^0)-f(\nabla u^0, u^0, x)=0,
\end{equation}
Note that from now on we shall always use $u^0$ to denote a solution of
a second order fully nonlinear PDE we intend to solve.
Equation \eqref{ma1} reduces to the classical Monge-Amp\`ere equation
\[
\mbox{det}(D^2 u^0)=f(x)
\]
if $f(\nabla u^0, u^0, x)=f(x) > 0$, and to Gauss curvature equation
\[
\mbox{det}(D^2 u^0)=K(1+|\nabla u^0|^2)^{\frac{n+2}2}
\]
if $f(\nabla u^0, u^0, x) =K(1+|\nabla u^0|^2)^{\frac{n+2}2}$.
Where the constant $K$ is a prescribed Gauss curvature.
Monge-Amp\`ere type equations are the best known second order fully
nonlinear PDEs, they arise in differential geometry and applications
such as mass transportation and meteorology. It is well-known that
Monge-Amp\`ere type equations are elliptic only in the set of convex
functions (cf. \cite{Gilbarg_Trudinger01,Gutierrez01}). So their
viscosity solutions are defined as convex functions in the
sense of Definition \ref{def2}.

The vanishing moment approximation \eqref{ePDE} to \eqref{ma1}
reads as:
\begin{equation}\label{eq5}
-\vepsi\Delta^2 u^\vepsi + \mbox{det}(D^2u^\vepsi)
=f(\nabla u^\vepsi, u^\vepsi, x) \qquad (\vepsi>0).
\end{equation}
For each fixed $\vepsi>0$, this is a quasilinear fourth order PDE
with Hessian type nonlinearity. It is complemented by boundary
conditions \eqref{bc2}$_1$ and \eqref{bc3}$_1$ (or \eqref{bc3}$_2$).

For the classical Monge-Amp\`ere equation, it can be shown that
\cite{feng06c,Neilan06}
\begin{itemize}
\item For each fixed $\vepsi>0$, problem \eqref{eq5}, \eqref{bc2}$_1$,
\eqref{bc3}$_1$ (or \eqref{bc3}$_2$) has a unique weak solution
in $W^{3,2}(\Ome)$ for $n=2, 3$.
\item $\mbox{det}(D^2 u^\vepsi) > 0$ and $\Del u^\vepsi > 0$
in $\Ome$ for $\vepsi>0$  when $n=2,3$.
\item $u^\vepsi$ is convex in $\Ome$ for $\vepsi >0$ when $n=2$.
\item The Dirichlet problem for the classical Monge-Amp\`ere equation
has a unique convex strong moment solution, which coincides with the
unique convex viscosity solution of the same problem when $n=2$.
\end{itemize}
The above results immediately imply that in the two dimensions
the vanishing moment method indeed works for the classical
Monge-Amp\`ere equation and the notion of moment solutions and
the notion of viscosity solutions are equivalent in this case.

\begin{remark}
Recall that when $n=2$, the Dirichlet problem \eqref{ma}--\eqref{mab}
has {\em at most} two solutions (see \cite{Courant_Hilbert89}). An amazing
numerical discovery to be
given in \S \ref{sec-5} is that if we restrict $\vepsi$ in \eqref{eq5}
to $\vepsi<0$, then $\underset{\vepsi\rightarrow 0^-}{\lim} u^\vepsi$
also exists and the limit is nothing but the other solution solution
of problem \eqref{ma}--\eqref{mab} which is {\em concave}
(see Figures \ref{fig1}--\ref{fig5})!
\end{remark}

\subsection{Pucci's equations} \label{sec-3.2a}
Pucci's extremal equations are referred to the following two families of 
fully nonlinear PDEs (cf. \cite{Gilbarg_Trudinger01,Caffarelli_Cabre95})
\begin{align}\label{pucci_1}
M_\alpha[u] &:= \alpha \Del u + (1-n\alpha) \lambda_n(D^2u^0) =f(x),\\
m_\alpha[u] &:= \alpha \Del u + (1-n\alpha) \lambda_1(D^2u^0) =f(x)
\label{pucci_2}
\end{align}
for $0< \alpha \leq \frac{1}{n}$. Where $\lambda_n(D^2u^0)$ 
and $\lambda_1(D^2u^0)$ denote the maximum and minimum eigenvalues 
of the Hessian matrix $D^2u^0$. In the $2$-d case,
the above equations can be rewritten in terms of $\Del u^0$ and 
$\mbox{det}(D^2 u^0)$ (cf. \cite{Dean_Glowinski06b}).

The vanishing moment approximations to \eqref{pucci_1} and 
\eqref{pucci_2} are defined as
\begin{align}\label{pucci_1a}
&-\vepsi \Del^2 u^\vepsi + \alpha \Del u^\vepsi 
+ (1-n\alpha) \lambda_n(D^2u^\vepsi) =f(x),\\
&-\vepsi \Del^2 u^\vepsi +\alpha \Del u^\vepsi 
+ (1-n\alpha) \lambda_1(D^2u^\vepsi) =f(x),
\label{pucci_2a}
\end{align}
which should be complemented by boundary conditions \eqref{bc2}$_1$
and \eqref{bc3}$_1$ (or \eqref{bc3}$_2$). 

\subsection{Infinite Laplace equation} \label{sec-3.2}
The infinite Laplace equation refers to the following {\em degenerate
quasilinear} PDE:
\begin{equation}\label{ip}
F(D^2u^0, Du^0, u^0,x):=\Del_\infty u^0 =0,
\end{equation}
where
\[
\Del_\infty u^0: = \langle D^2 u^0 \nabla u^0,\nabla u^0\rangle
= D^2 u^0 \nabla u^0\cdot \nabla u^0.
\]
$\Del_\infty u^0$ can be regarded as the limit of the $p$-Laplacian
$\Del_p u^0:=\Div(|\nabla u^0|^{p-2}\nabla u^0)$ as $p\rightarrow \infty$,
it also can be derived as the Euler-Lagrange equation of the
$L^\infty$ functional
\[
I(v):= \mbox{ess sup}_{x\in \Ome} |\nabla v(x)|,
\]
whose minimizers are often called ``absolute minimizers"
\cite{Aronsson_Crandall_Juutinen04}.
Besides its mathematical appeals, the infinite Laplace equation
also arises from image processing, geography, and geology
applications \cite{Aronsson_Evans_Wu96,Caselles_Morel_Sbert98}.
Although the infinite Laplace equation is only a degenerate
quasilinear PDE, not a fully nonlinear PDE, it is very
difficult to solve numerically. This is because the infinite Laplace
equation does not have classical solutions in general
\cite{Aronsson_Crandall_Juutinen04}, and since it is not in
divergence form, its weak solutions are defined and understood
in the viscosity sense. We refer to \cite{Oberman05} for 
recent developments on finite difference approximations of
the infinite Laplace equation.

Here we propose the following vanishing moment approximation for \eqref{ip}:
\begin{equation}\label{ipe}
-\vepsi \Del^2 u^\vepsi + \Del_\infty u^\vepsi =0,
\end{equation}
which is complemented by boundary conditions \eqref{bc2}$_1$
and \eqref{bc3}$_1$ (or \eqref{bc3}$_2$).
In \S \ref{sec-5}, we shall present numerical results which
show that the vanishing moment approximation exactly converges
to the unique viscosity solution of the Dirichlet
problem for \eqref{ip}. This is another example which shows that
the notion of moment solutions and the notion of viscosity
solutions coincide.

\begin{remark}
It is easy to see that the above vanishing moment method also applies to
the $p$-Laplacian equation $-\Del_p u^0=f$ for $1\leq p< \infty$.
\end{remark}

\subsection{Second order fully nonlinear parabolic PDEs} \label{sec-3.3}
We first like to note that there are several different versions of
legitimate parabolic generalizations to elliptic PDE \eqref{eq2}
(cf. \cite{Lieberman96,Wang92a}). In this paper, we shall
only consider the following widely studied (and it turns out
to be the ``easiest") class of second order fully nonlinear parabolic PDEs:
\begin{equation}\label{parabolic1}
F(D^2u^0, \nabla u^0, u^0, x, t)-u^0_t=0,
\end{equation}
assuming that $F(D^2u^0, \nabla u^0, u^0, x,t)$ is elliptic.
Clearly, this is the most natural parabolic generalization
to equation \eqref{eq2}. For example, the corresponding
parabolic Monge-Amper\'e type equation reads as
\begin{equation}\label{pma1}
\mbox{det}(D^2u^0)-u^0_t=f(\nabla u^0,u^0,x,t) \geq 0.
\end{equation}
In the past two decades the viscosity solution theory
has been well developed for equations \eqref{parabolic1}
and \eqref{pma1}, see
\cite{Lieberman96,Wang92a,Gutierrez_Huang01}.
On the other hand, numerical approximation to these fully
nonlinear parabolic PDEs is a completely untouched area. To
the best our knowledge, no numerical result (in fact, no attempt)
is known in the literature.

Similarly, we can define the vanishing moment method and
the notion of moment solutions for initial and initial-boundary
value problems for \eqref{parabolic1}, and then ask the same
questions as we did in \S \ref{sec-2.2}. We leave this as an
exercise to interested readers and refer to \cite{Neilan06}
for a detailed exposition.

Following the derivation of \S \ref{sec-2.2}, we propose
the following vanishing moment approximations to \eqref{parabolic1}
and \eqref{pma1}, respectively,
\begin{eqnarray}\label{parabolic2}
F(D^2u^\vepsi, \nabla u^\vepsi, u^\vepsi, x, t)
- \vepsi\Del^2 u^\vepsi - u^\vepsi_t &=& 0,\\
\mbox{det}(D^2u^\vepsi) - \vepsi\Del^2 u^\vepsi - u^\vepsi_t &=&
f(\nabla u^\vepsi, u^\vepsi, x, t),
\label{pma2}
\end{eqnarray}
each of the above equations is a fourth order quasilinear
parabolic PDEs. 

\section{Numerical experiments} \label{sec-5}

In this section, we shall present a number of numerical 
experiment results obtained by using the vanishing moment
method together with the numerical methods proposed in
\S\ref{sec-4}. Both $2$-d and $3$-d tests will be 
presented. All the $3$-d tests are obtained by a 
Hermann-Miyoshi type mixed finite element method,
while the $2$-d tests are computed by using 
both the Argyris (plate) finite element method
and the Hermann-Miyoshi mixed finite element method. 

\subsection{Two-dimensional numerical experiments}\label{sec-5.1}

The numerical solutions of the first seven tests are computed 
using the Argyris finite element method. 

{\bf Test 1:} In this test we solve the Monge-Amp\'ere problem
\eqref{ma}--\eqref{mab} on the unit square $\Ome=(0,1)^2$
with the following data:
\[
f(x,y)\equiv 1, \qquad g(x,y)\equiv 0.
\]
We remark that problem \eqref{ma}--\eqref{mab} has 
a unique convex viscosity solution but does not 
have a classical solution (cf. \cite{Gutierrez01,Dean_Glowinski06b}). 

\begin{figure}[htb]
\centerline{
\includegraphics[angle=0,width=7cm,height=5cm]{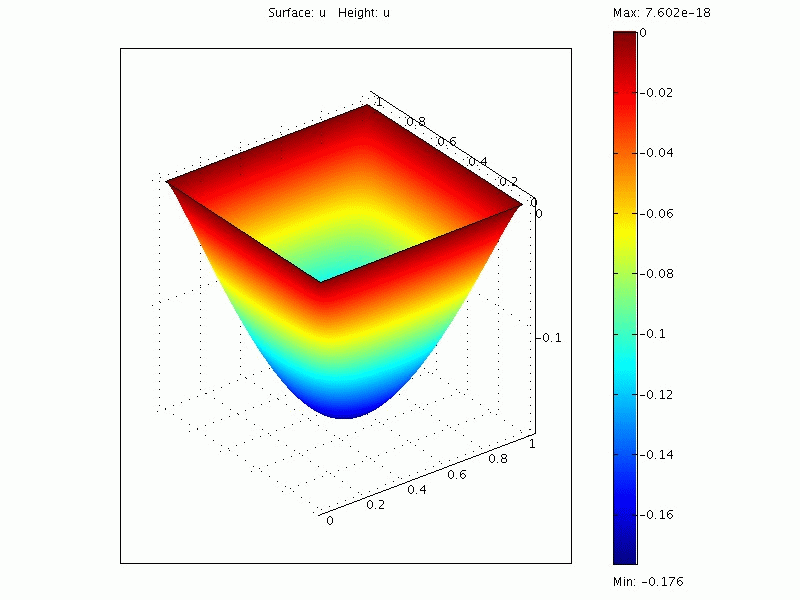}
\includegraphics[angle=0,width=7cm,height=5cm]{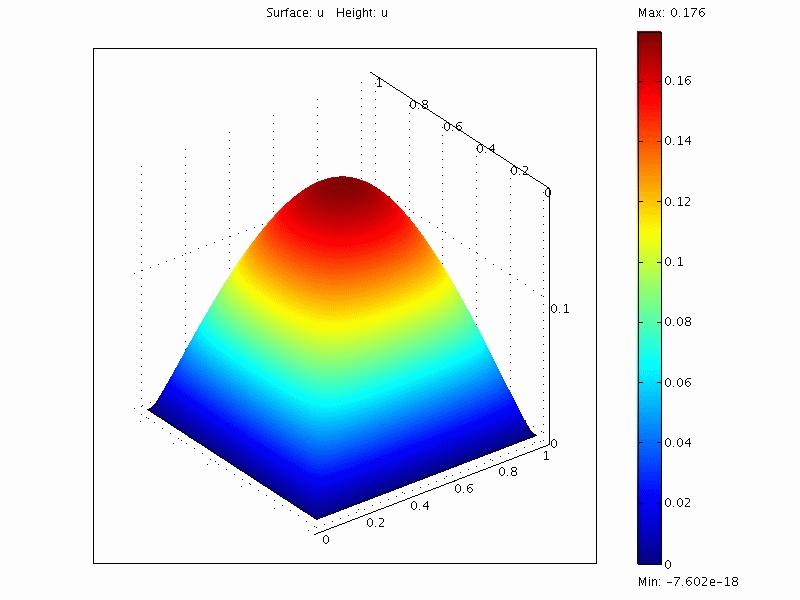}
}
\caption{{\small Computed (moment) solutions of Test 1: 
Graph on left corresponds to $\vepsi=10^{-3}$ and graph on 
right corresponds to $\vepsi=-10^{-3}$.}}\label{fig1}
\end{figure}

\begin{figure}[htb]
\centerline{
\includegraphics[angle=0,width=6.5cm,height=5cm]{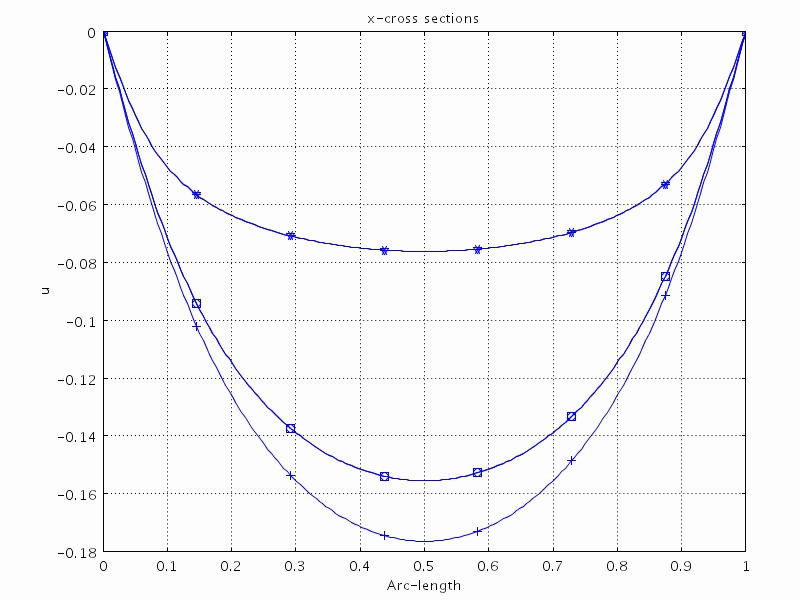}
\includegraphics[angle=0,width=6.5cm,height=5cm]{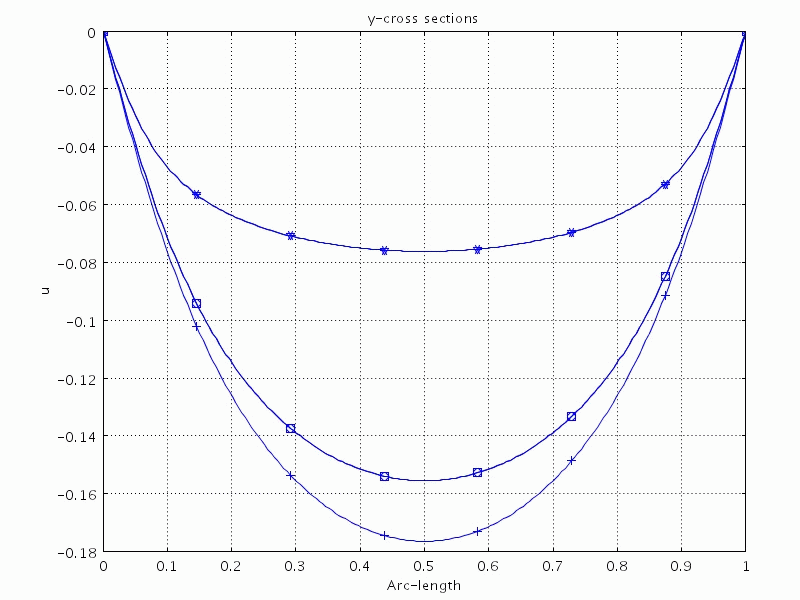}
}
\caption{{\small x-cross sections (left figure) of the left graph in Fig. \ref{fig1} at
$x=0.1, 0.3, 0.5, 0.7, 0.9$ (indicated respectively by asterisk, circle, plus sign,
squre, and triangle); y-cross sections (right figure) of the left graph in Fig. \ref{fig1} 
at $y=0.1, 0.3, 0.5, 0.7, 0.9$ (indicated respectively by asterisk, circle, plus sign,
squre, and triangle)}}\label{fig1.1}
\end{figure}

Recall that the vanishing moment approximation of \eqref{ma}--\eqref{mab}
is problem \eqref{eq5}, \eqref{bc2}$_1$, \eqref{bc3}$_1$ with the above
$f$ and $g$. We discretize problem \eqref{eq5}, \eqref{bc2}$_1$, 
\eqref{bc3}$_1$ using the Argyris plate element as described in 
\S \ref{sec-4.1}. Figure \ref{fig1} displays the computed (moment)
solutions using $\vepsi=10^{-3}$ (left graph) and $\vepsi= -10^{-3}$ 
(right graph). Clearly, the vanishing moment approximations 
correctly capture the convex viscosity solution (left graph)
and the concave viscosity solution (right graph). Hence, the 
moment solutions coincide with the viscosity solutions (see
\cite{feng06c} for a rigorous proof).

To have a better view of the convexity of the computed solution, we also
plot selected cross sections of the left figure in Figure \ref{fig1}.  The cross
sections clearly show that the computed solution is a convex function. In particular,
there is no visible smear at the boundary.

{\bf Test 2:} The only difference between this test and Test 1 is that
the datum functions are now chosen as
\[
f(x,y)=(1+(x^2+y^2))e^{(x^2+y^2)},\quad
g(x,y)=\left\{
       \begin{array}{ll}
        e^{y^2/2} &\quad\mbox{if } x=0, \\
        e^{x^2/2} &\quad\mbox{if } y=0,\\
        e^{(1+x^2)/2} &\quad\mbox{if } y=1,\\
        e^{(1+y^2)/2} &\quad\mbox{if } x=1,
       \end{array} \right.
\]
so that $u^0(x,y)= \frac12 e^{(x^2+y^2)}$ is an exact solution 
of problem \eqref{ma}--\eqref{mab}. Clearly, $u^0$ is a convex function, 
hence $u^0$ must be the unique convex viscosity solution of 
problem \eqref{ma}--\eqref{mab} (cf. \cite{Gutierrez01}).

\begin{figure}[htb]
\centerline{
\includegraphics[angle=0,width=7cm,height=5cm]{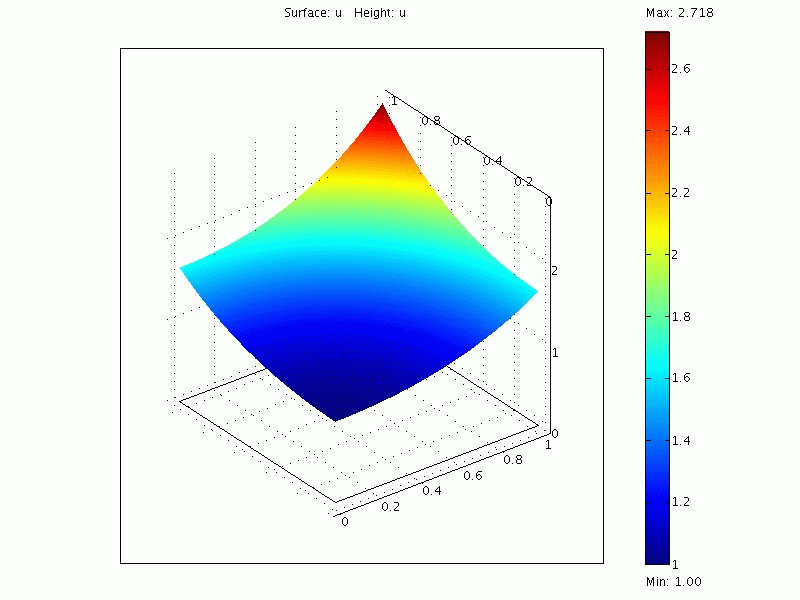}
\includegraphics[angle=0,width=7cm,height=5cm]{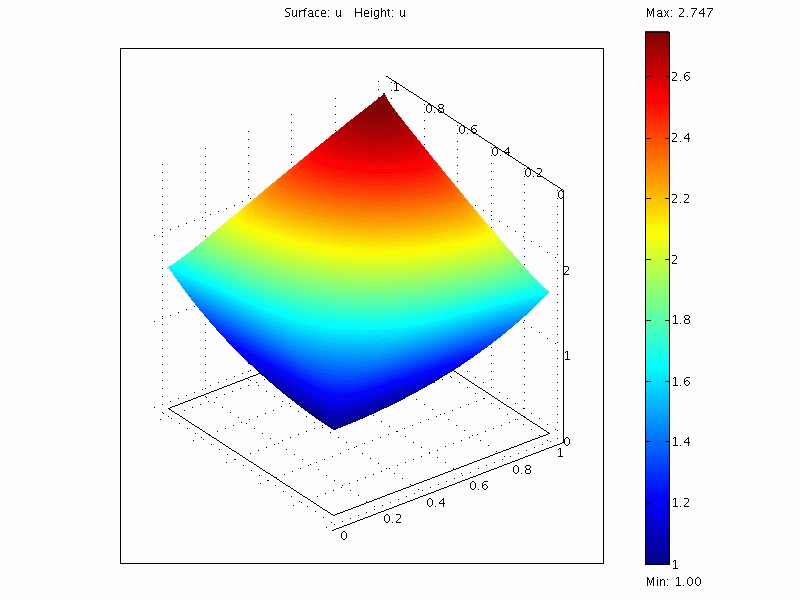}
}
\caption{{\small Computed (moment) solutions of Test 2:
Graph on left corresponds to $\vepsi=10^{-3}$ and graph on 
right corresponds to $\vepsi=-10^{-3}$. }}\label{fig2}
\end{figure}

Figure \ref{fig2} shows the computed (moment)
solutions using $\vepsi=10^{-3}$ (left graph) and $\vepsi= -10^{-3}$
(right graph). Again, the vanishing moment approximations
correctly capture the convex viscosity solution $u^0$ (left graph)
and the concave viscosity solution (right graph), hence the
moment solutions coincide with the viscosity solutions (see
\cite{feng06c} for a rigorous proof).

{\bf Test 3:} Similar to Test 2, the only difference between this 
test and Test 1 is that the datum functions are now chosen as
\[
f(x,y)=\frac{1}{x^2+y^2},\qquad
g(x,y)=\left\{
  \begin{array}{ll}
  \frac{2\sqrt{2}}{3} y^{\frac32} &\quad\mbox{if } x=0,\\
  \frac{2\sqrt{2}}{3} x^{\frac32} &\quad\mbox{if } y=0,\\
  \frac{2\sqrt{2}}{3} (1+x^2)^{\frac34} &\quad\mbox{if } y=1,\\
  \frac{2\sqrt{2}}{3} (1+y^2)^{\frac34} &\quad\mbox{if } x=1,
  \end{array} \right.
\]
so that $u^0(x,y)= \frac{2\sqrt{2}}{3} (x^2+y^2)^{\frac34}$ 
is the unique convex viscosity solution of problem \eqref{ma}--\eqref{mab}. 

\begin{figure}[htb]
\centerline{
\includegraphics[angle=0,width=7cm,height=5cm]{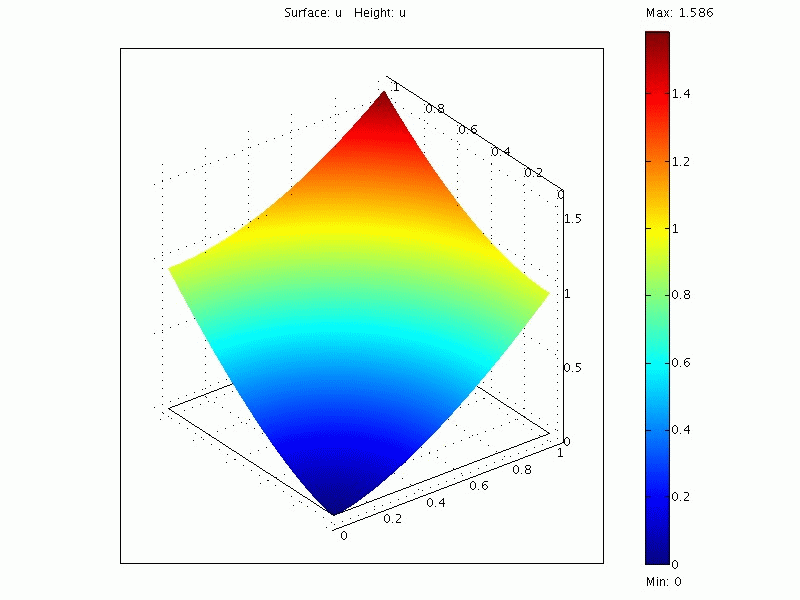}
\includegraphics[angle=0,width=7cm,height=5cm]{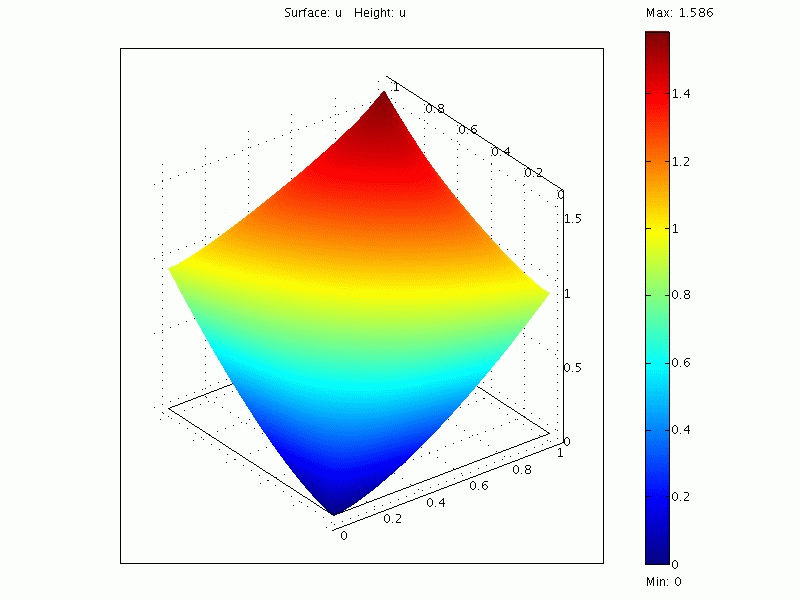}
}
\caption{{\small Computed (moment) solutions of Test 3:
Graph on left corresponds to $\vepsi=10^{-3}$ and graph on 
right corresponds to $\vepsi=-10^{-3}$. }}\label{fig3}
\end{figure}

Figure \ref{fig3} displays the computed (moment)
solutions using $\vepsi=10^{-3}$ (left graph) and $\vepsi= -10^{-3}$
(right graph). As expected, the vanishing moment approximations
correctly capture the convex viscosity solution $u^0$ (left graph)
and the concave viscosity solution (right graph), hence the
moment solutions coincide with the viscosity solutions (see
\cite{feng06c} for a rigorous proof).

{\bf Test 4:} Again, the only difference between this
test and Test 1 is that the datum functions are now chosen as
\[
f(x,y) = (1-x-y)^2 \qquad g\equiv 0 .
\]
On the other hand, mathematically there is a significant difference 
between these two test problems. Note that $f(x,y)=0$ on the line
$x+y=1$ in the domain $\Ome=(0,1)^2$. Hence, problem \eqref{ma}--\eqref{mab}
is known as a {\em degenerate} Monge-Amp\'ere problem (cf. \cite{Gutierrez01}).

\begin{figure}[htb]
\centerline{
\includegraphics[angle=0,width=7cm,height=5cm]{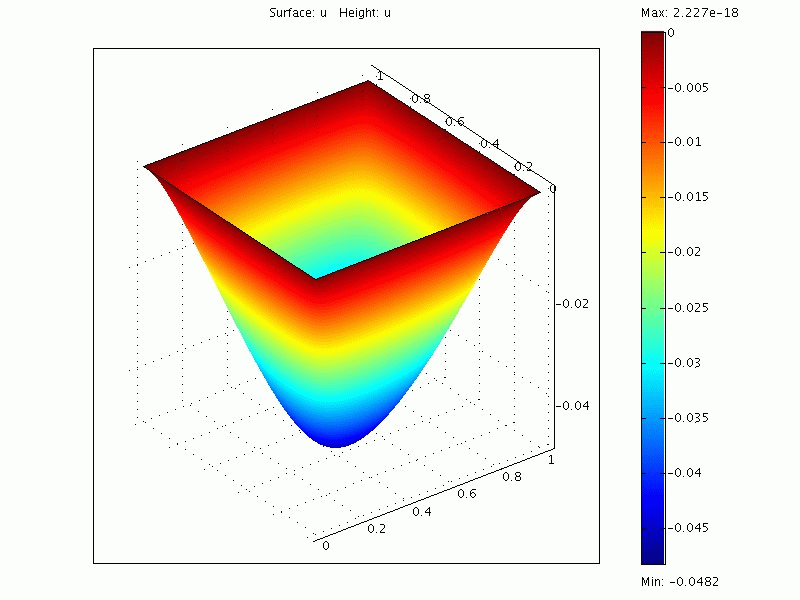}
\includegraphics[angle=0,width=7cm,height=5cm]{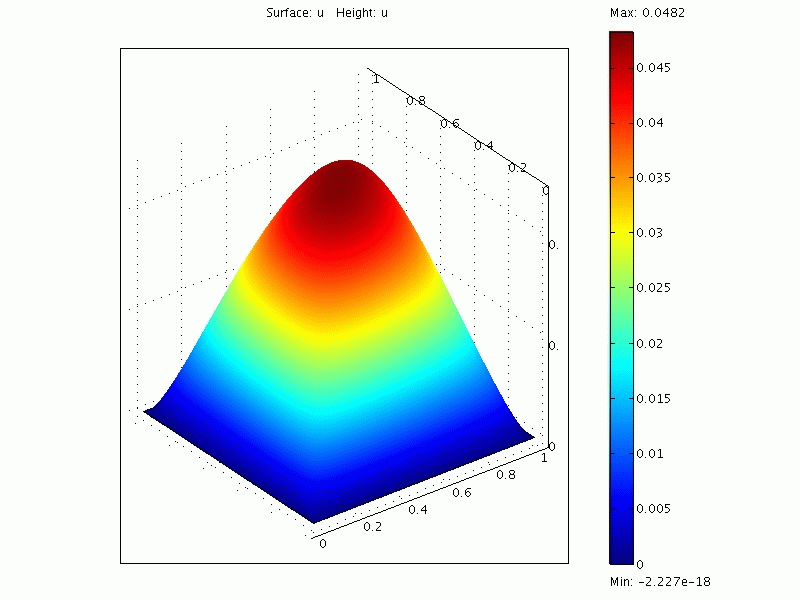}
}
\caption{{\small Computed (moment) solutions of Test 4:
Graph on left corresponds to $\vepsi=10^{-3}$ and graph on 
right corresponds to $\vepsi=-10^{-3}$. }}\label{fig4}
\end{figure}

Figure \ref{fig4} displays the computed (moment)
solutions using $\vepsi=10^{-3}$ (left graph) and $\vepsi= -10^{-3}$
(right graph). Once again, the vanishing moment approximations
correctly capture the convex viscosity solution (left graph)
and the concave viscosity solution (right graph), hence the
moment solutions coincide with the viscosity solutions (see
\cite{feng06c} for a rigorous proof). In addition, our numerical
result shows that the vanishing moment method is robust with
respect to the degeneracy of the underlying PDE.

{\bf Test 5:} Once again, the only difference between this
test and Test 1 is that the datum functions are now chosen as
\[
f(x,y)= x^2-y^2 \qquad  g\equiv 0 .
\]
Mathematically, the difference between this test problem and
Test 1 is even more dramatic because not only $f(x,y)=0$ on the 
line $x-y=0$ {\em but also $f$ changes sign} (hence the PDE changes 
type) in $\Ome$. To the best of our knowledge, there is no 
viscosity solution theory for this type Monge-Ampere problems 
in the literature. However, the vanishing moment method
seem works well for this problem. Our numerical results
indicate existence of  both convex and concave moment solutions.

\begin{figure}[htb]
\centerline{
\includegraphics[angle=0,width=7cm,height=5cm]{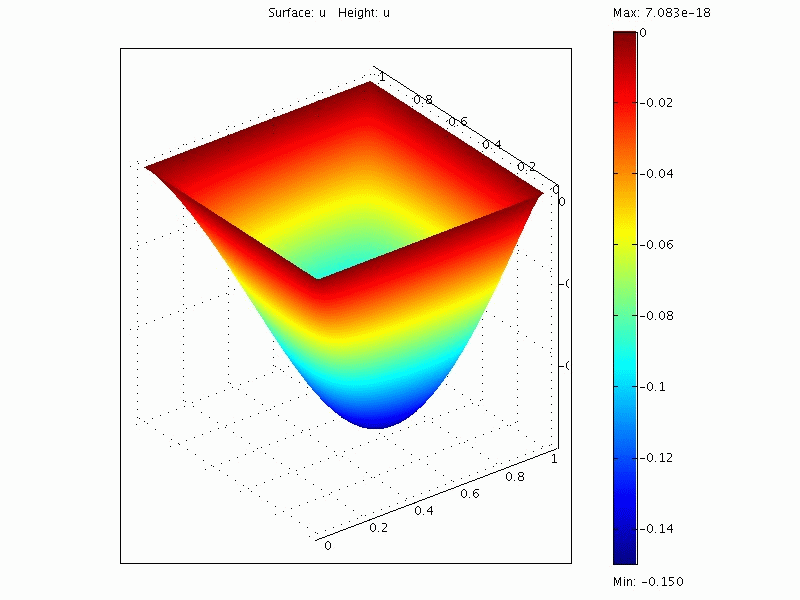}
\includegraphics[angle=0,width=7cm,height=5cm]{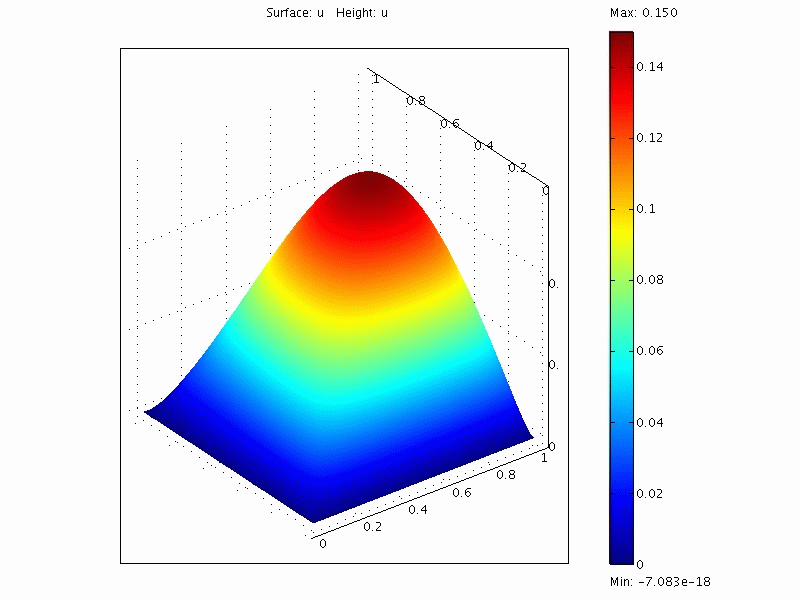}
}
\caption{{\small Computed (moment) solutions of Test 5:
Graph on left corresponds to $\vepsi=10^{-3}$ and graph on 
right corresponds to $\vepsi=-10^{-3}$. }}\label{fig5}
\end{figure}

Figure \ref{fig5} displays the computed convex (moment) 
solution using $\vepsi=10^{-3}$ (left graph) and the computed 
concave (moment) solution using $\vepsi= -10^{-3}$ (right graph).  

Again, to have a better view of the convexity of the computed solution, we also
plot selected cross sections of the left figure in Figure \ref{fig5}.  The cross
sections clearly show that the computed solution is a convex function. In particular,
there is no visible smear at the boundary. 
\begin{figure}[htb]
\centerline{
\includegraphics[angle=0,width=6.5cm,height=5cm]{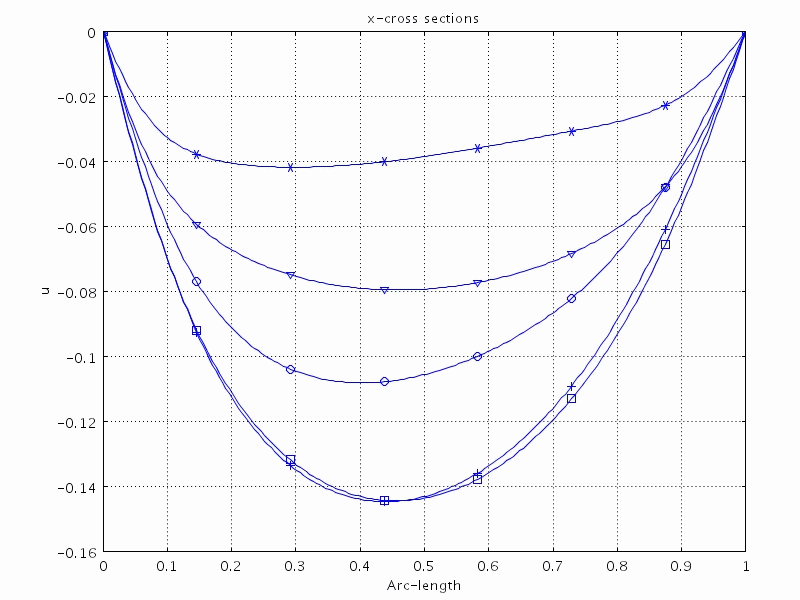}
\includegraphics[angle=0,width=6.5cm,height=5cm]{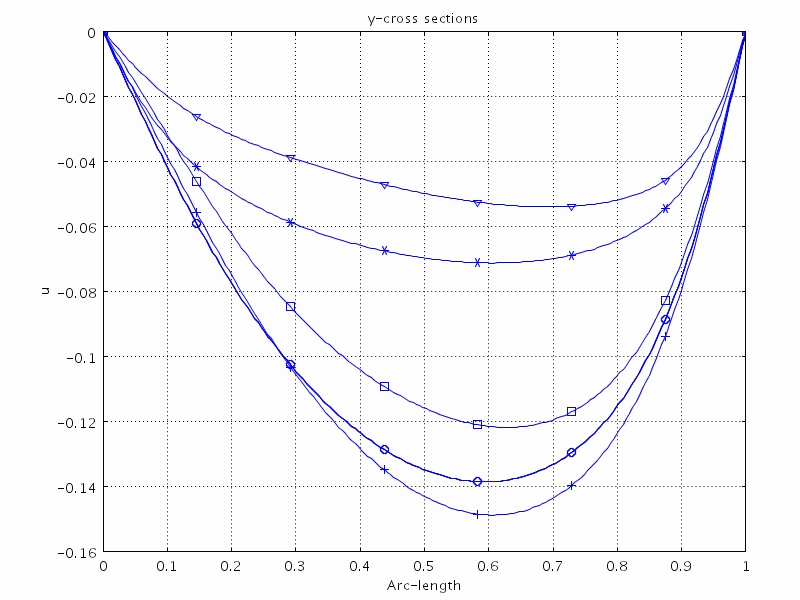}
}
\caption{{\small x-cross sections (left figure) of the left graph in Fig. \ref{fig5} at 
$x=0.1, 0.3, 0.5, 0.7, 0.9$ (indicated respectively by asterisk, circle, plus sign,
squre, and triangle); y-cross sections (right figure) of the left graph in Fig. \ref{fig5} 
at $y=0.1, 0.3, 0.5, 0.7, 0.9$ (indicated respectively by asterisk, circle, plus sign,
squre, and triangle)}}\label{fig5.1}
\end{figure}

{\bf Test 6:} In this test we solve the following Gauss curvature 
(or $\mathcal{K}$-surface) equation (cf. 
\cite{Guan95,Guan_Spruck04,Guan_Guan02})
\begin{alignat}{2} \label{gauss1}
\mbox{det} (D^2 u^0) &= K (1+|\nab u^0|^2)^2 &&\qquad\mbox{in } 
\Ome:=(-0.57,0.57)^2, \\
u^0&=x^2+y^2-1 &&\qquad\mbox{on } \p\Ome, \label{gauss2}
\end{alignat}
where $K>0$ is a given constant Gauss curvature.
Note that the above problem is a special case of 
problem \eqref{ma1},\eqref{bc1} with 
$f(\nab u^0,u^0,x,y)= K(1+|\nab u^0|^2)^{\frac{n+2}2}$, $n=2$, and
$g(x,y)=x^2+y^2-1$. 

It was proved by Guan \cite{Guan95} that 
there exists $K^*>0$ such that for each $K\in [0, K^*)$ 
problem \eqref{gauss1}--\eqref{gauss2} (with more general
Dirichlet data) has a unique convex viscosity solution.  
Theoretically, it is very difficult to give an accurate 
estimate for the curvature upper bound $K^*$. However,
hand, this offers an ideal opportunity for numerical analysts 
to help and to contribute. It turns out that the 
vanishing moment method proposed in this paper 
works very well for such a problem, hence it might 
provide a useful tool and answer to the challenge.  

\begin{figure}[htb]
\centerline{
\includegraphics[angle=0,width=7cm,height=5cm]{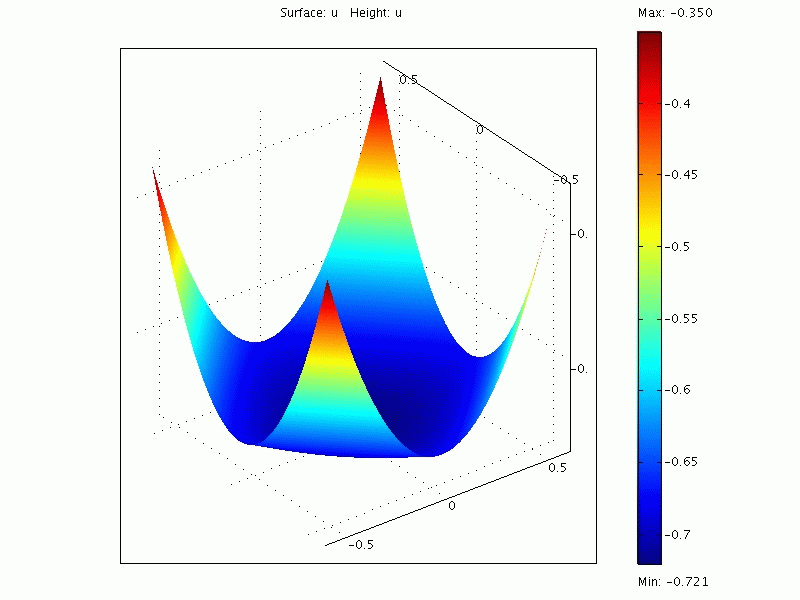}
\includegraphics[angle=0,width=7cm,height=5cm]{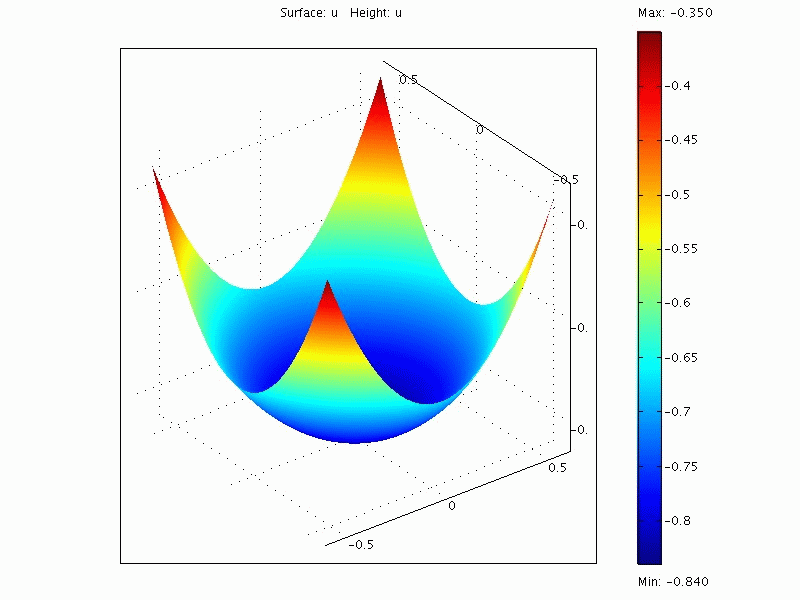}
}
\centerline{
\includegraphics[angle=0,width=7cm,height=5cm]{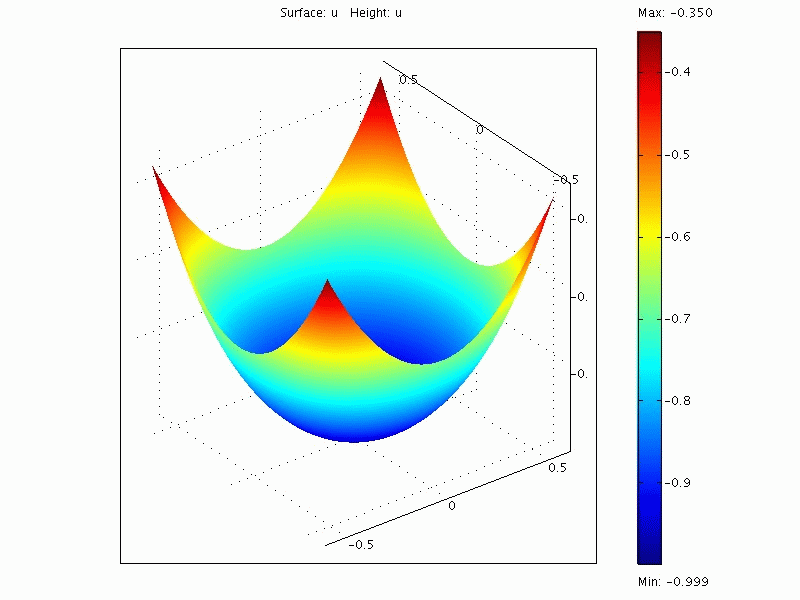}
\includegraphics[angle=0,width=7cm,height=5cm]{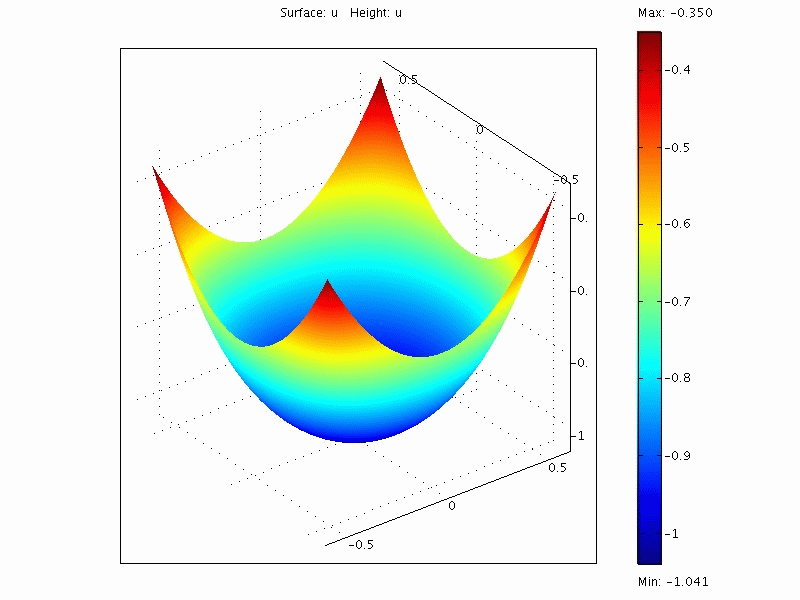}
}
\caption{{\small  Computed (moment) solutions of Test 6: $\vepsi=10^{-3}$ and
$K=0.1, 1, 2, 2.1$. Graphs are arranged row-wise.}}\label{fig6}
\end{figure}

Since we are only interested in convex solutions of the Gauss curvature
equation, so we restrict $\vepsi>0$ in \eqref{eq5}. Figure \ref{fig6} 
displays the computed convex (moment) solution using $\vepsi=10^{-3}$ 
and $K=0.1, 1, 2, 2.1$, respectively.
We note that our computer code stops producing a convergent numerical 
solution for $K=2.2$.  Hence we conjecture that $K^*\approx 2.1$ 
for the above test problem. 

{\bf Test 7:} In this test, we solve problem \eqref{ip},
\eqref{bc1} over the domain $\Ome=(-\frac12,\frac12)^2$ with the 
following boundary datum function
\[
g(x,y)=\left\{
  \begin{array}{ll}
  \Bigl(\frac12\Bigr)^{\frac43}- y^{\frac43} &\quad\mbox{if } x=-\frac12,\\
  x^{\frac43}-\Bigl(\frac12\Bigr)^{\frac43} &\quad\mbox{if } y=-\frac12,\\
  \Bigl(\frac12\Bigr)^{\frac43}- y^{\frac43} &\quad\mbox{if } x=\frac12,\\
  x^{\frac43}-\Bigl(\frac12\Bigr)^{\frac43}  &\quad\mbox{if } y=\frac12,
  \end{array} \right.
\]
so that $u^0(x,y)=x^{\frac43}-y^{\frac43}$ is the unique viscosity solution
(cf. \cite{Aronsson_Crandall_Juutinen04}).

We remark that this is an
important example in the theory of absolutely minimizing
functions since $u^0$ is the least regular absolutely minimizing
function known in the case of the Euclidean norm 
(cf. \cite{Aronsson_Crandall_Juutinen04} and the references therein). 
It is easy to check that $u^0$ is a H\"older continuous
function with exponent $\frac13$. However, it is not
twice differentiable on the axes.
We also note that (see \S \ref{sec-3.2}) the infinite Laplace
equation \eqref{ip} is only a second order (degenerate)
quasilinear (instead fully nonlinear) PDE,  however, the
complicate and nondivergence structure makes  
the infinite Laplace equation very difficult to 
analyze theoretically and to compute numerically.

\begin{figure}[t] 
\centerline{
\includegraphics[angle=0,width=7cm,height=5cm]{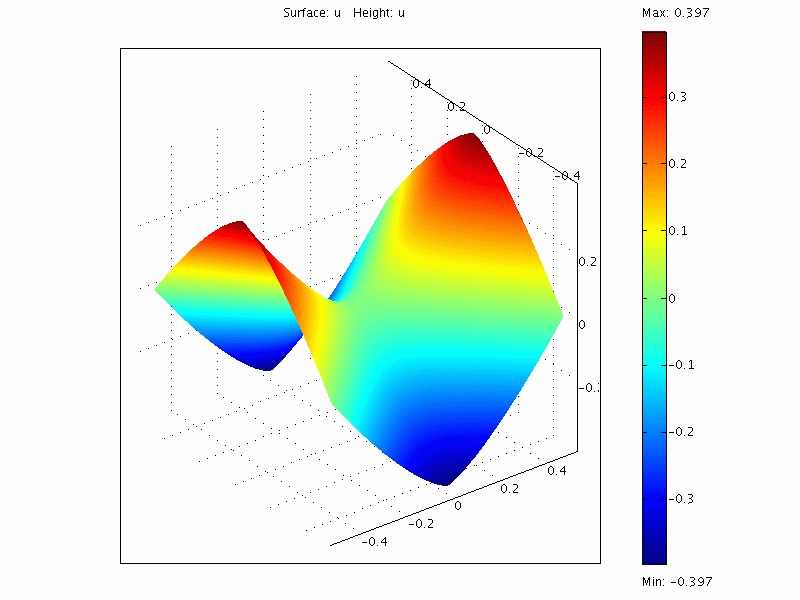}
\includegraphics[angle=0,width=7cm,height=5cm]{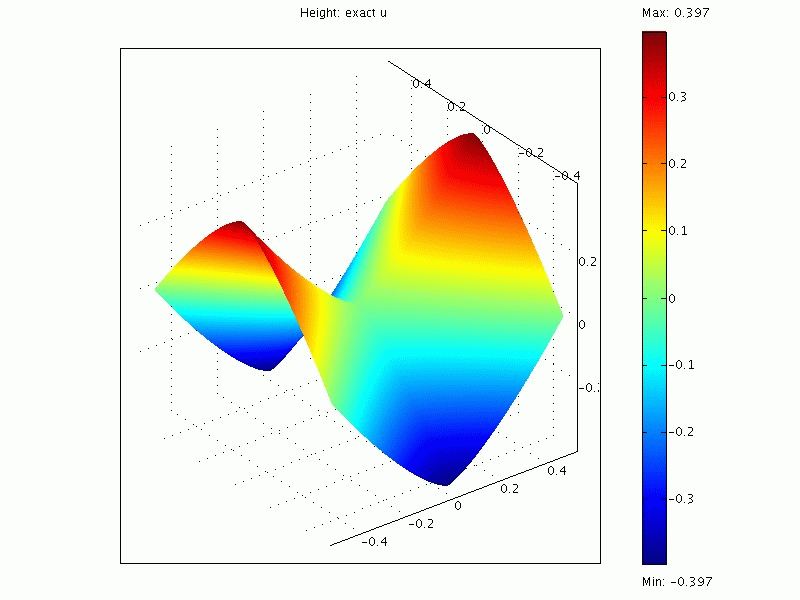}
}
\caption{{\small  Computed moment solution (left graph) and the exact 
viscosity solution (right graph) of Test 7. $\vepsi=10^{-3}$.}}\label{fig7}
\end{figure}

Figure \ref{fig7} displays the computed (moment)
solution using $\vepsi=10^{-3}$ (left graph) and the exact solution $u^0$
(right graph). Once again, the vanishing moment approximation
correctly captures the viscosity solution $u^0$, hence, the
moment solution coincides with the viscosity solution (see
\cite{feng06c} for a rigorous proof). 

The numerical solutions of the next two tests are obtained by using the
Hermann-Miyoshi mixed finite element method with piecewise quadratic 
shape functions.

\begin{figure}[htb]
\centerline{
\includegraphics[angle=0,width=7cm,height=5cm]{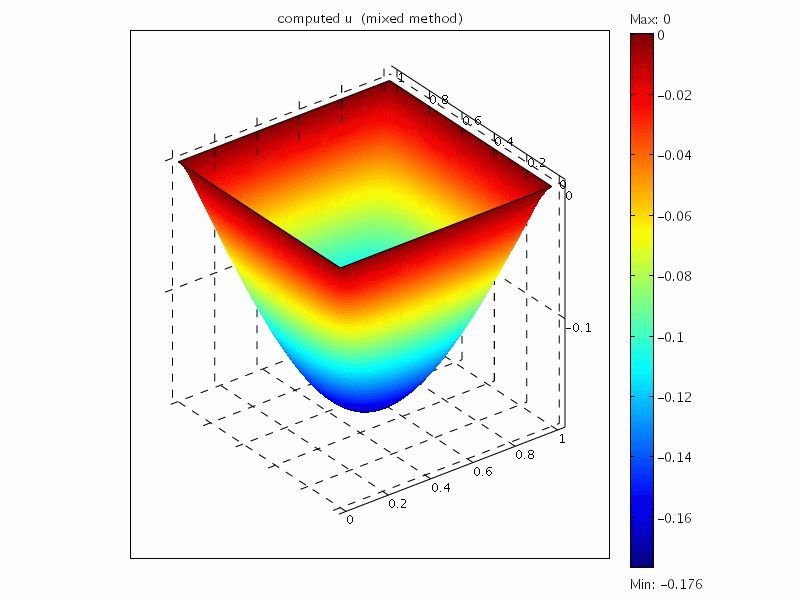}
\includegraphics[angle=0,width=7cm,height=5cm]{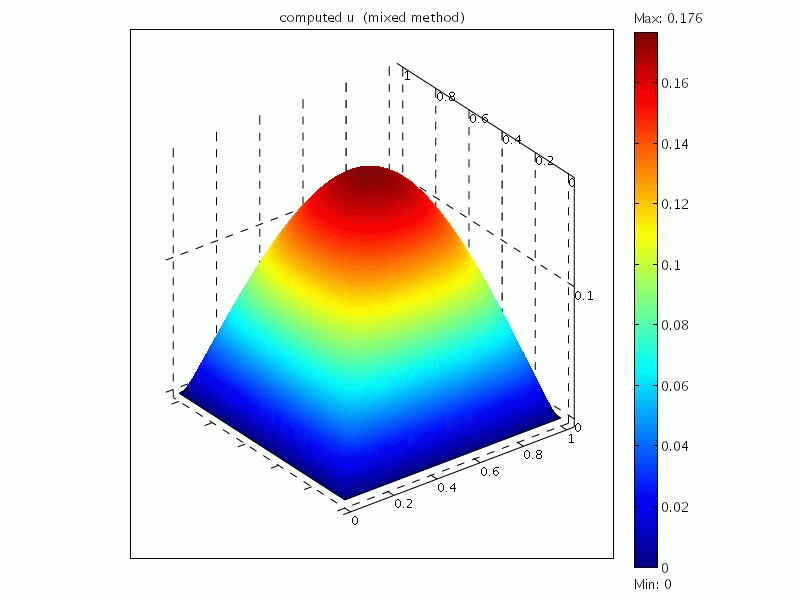}
}
\caption{{\small Computed (moment) solutions of Test 8:
Graph on left is the computed $u^\vepsi_h$ with $\vepsi=10^{-3}$ 
and graph on right is the computed $u^\vepsi_h$ with 
$\vepsi=-10^{-3}$.}}\label{fig8}
\end{figure}

\smallskip
{\bf Test 8:}  This test is a re-run of Test 1 but using the quadratic 
Hermann-Miyoshi mixed finite element method. Figure \ref{fig8} 
is the counterpart of Figure \ref{fig1}.  
We remark that the mixed method also produces an approximation
to the Hessian matrix $D^2u^\vepsi$, which is not shown here. Clearly,
the numerical results of Test 1 and Test 8 have the same accuracy.
However, it should be noted that the mixed method 
runs about $20$ times faster than the Argyris method on this test 
problem.

\smallskip
{\bf Test 9:} This test solves, using the quadratic Hermann-Miyoshi
mixed finite method, the Monge-Amp\'ere  problem \eqref{ma}--\eqref{mab} 
on the unit square $\Ome=(0,1)^2$ with the following data:
\[
f(x,y)=\frac{4}{(4-x^2-y^2)^2}, \qquad g(x,y)= \sqrt{4-x^2-y^2}
\]
so that $u^0=\sqrt{4-x^2-y^2}$ is an exact (convex) solution. 
We note that problem \eqref{ma}--\eqref{mab} 
has exact two solutions, one is convex and the other is concave
(cf.\cite{Courant_Hilbert89}).

\begin{figure}[t]
\centerline{
\includegraphics[angle=0,width=7cm,height=5cm]{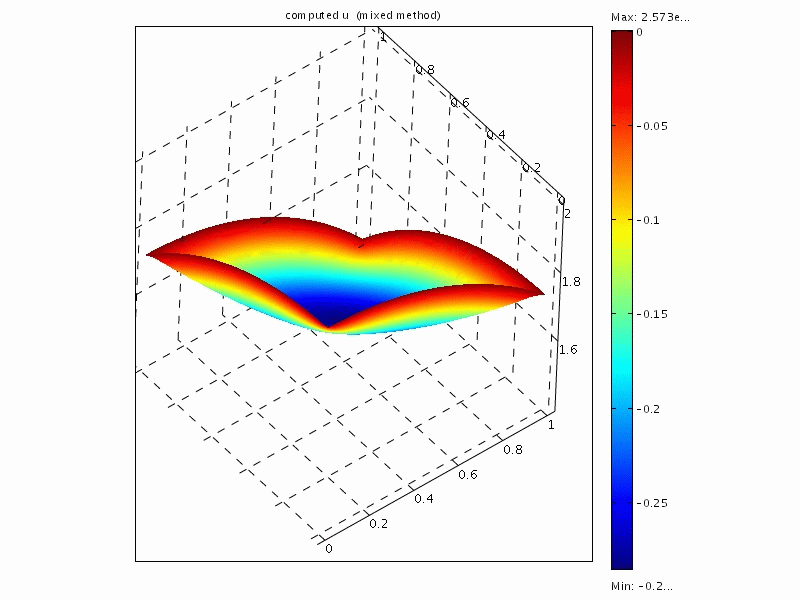}
\includegraphics[angle=0,width=7cm,height=5cm]{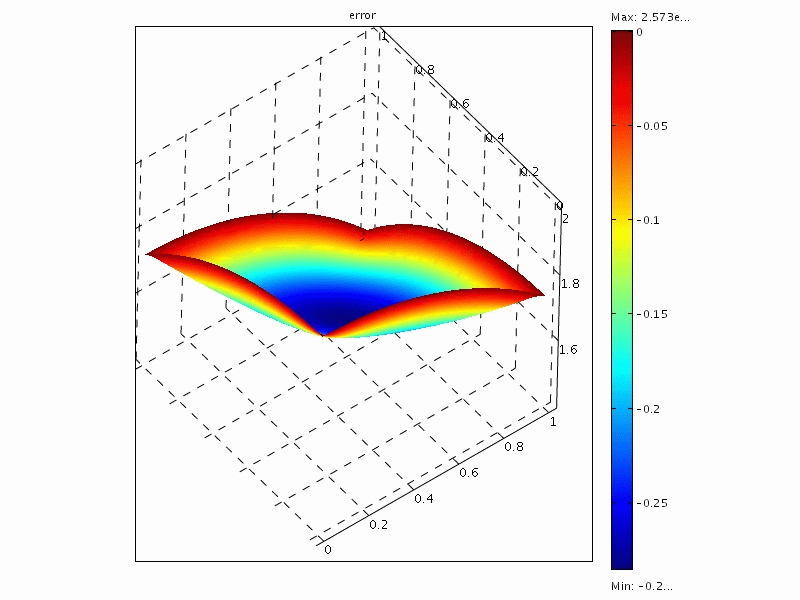}
}
\caption{{\small Computed (moment) solutions of Test 9 by
the quadratic Hermann-Miyoshi mixed method. $\vepsi=10^{-3}$.}}\label{fig9}
\end{figure}

Figure \ref{fig9} displays the computed (moment)
solution $u^\vepsi_h$ using $\vepsi=10^{-3}$ (left graph) and its
error (right graph). As expected, the vanishing moment approximation
correctly captures the convex viscosity solution $u^0$. Hence the
moment solution coincides with the viscosity solution (see
\cite{feng06c} for a rigorous proof). 
Again, we remark that the mixed method also gives
an approximation to the Hessian matrix $D^2u^\vepsi$, 
which is not shown here, and the mixed method
runs about $20$ times faster than the Argyris method 
for solving this test problem.

\subsection{Three-dimensional numerical experiments}\label{sec-5.2}
In this subsection we present two numerical tests on computing
moment (and viscosity) solutions of the Monge-Amp\'ere  
problem \eqref{ma}--\eqref{mab} in the unit cube $\Ome=(0,1)^3$. 
Numerical approximations of fully nonlinear PDEs in $3$-d is 
known to be very difficult. To the best of our knowledge,  
no $3$-d numerical results are given in the literature for 
the Monge-Amp\'ere type fully nonlinear PDEs.  

\smallskip
{\bf Test 10:} Consider the Monge-Amp\'ere problem \eqref{ma}--\eqref{mab}
on the unit cube $\Ome=(0,1)^3$ with the following data:
\[
f(x,y,z)=(1+x^2+y^2+x^2) \exp\Big(\frac{x^2+y^2+x^2}{2}\Bigr), 
\, g(x,y,z)= \exp\Big(\frac{x^2+y^2+x^2}{2}\Bigr).
\]
It is easy to verify that $u^0=\exp\big(\frac{x^2+y^2+x^2}{2}\bigr)$ 
is a unique exact (convex) solution. We compute this solution
using the vanishing moment method combined with a generalized 
Hermann-Miyoshi type mixed finite element method using linear 
shape functions.

\begin{figure}[htb]
\centerline{
\includegraphics[angle=0,width=7cm,height=5cm]{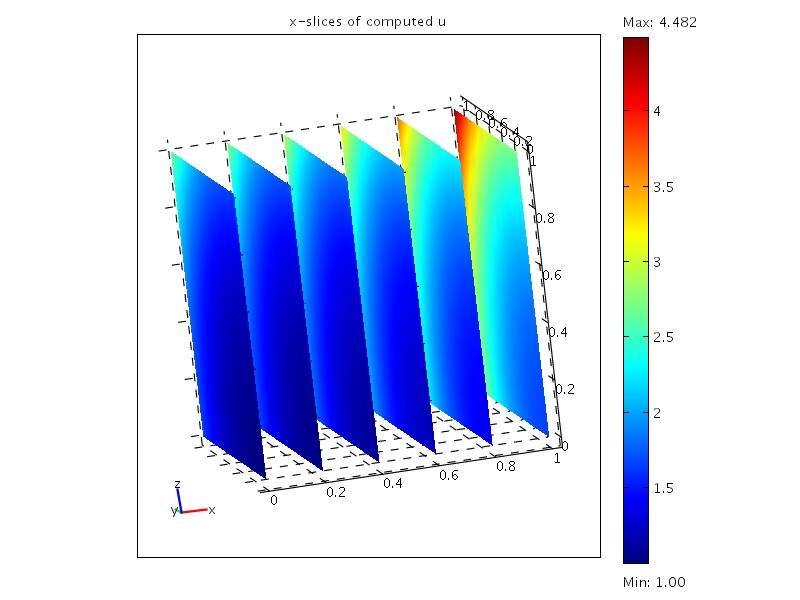}
\includegraphics[angle=0,width=7cm,height=5cm]{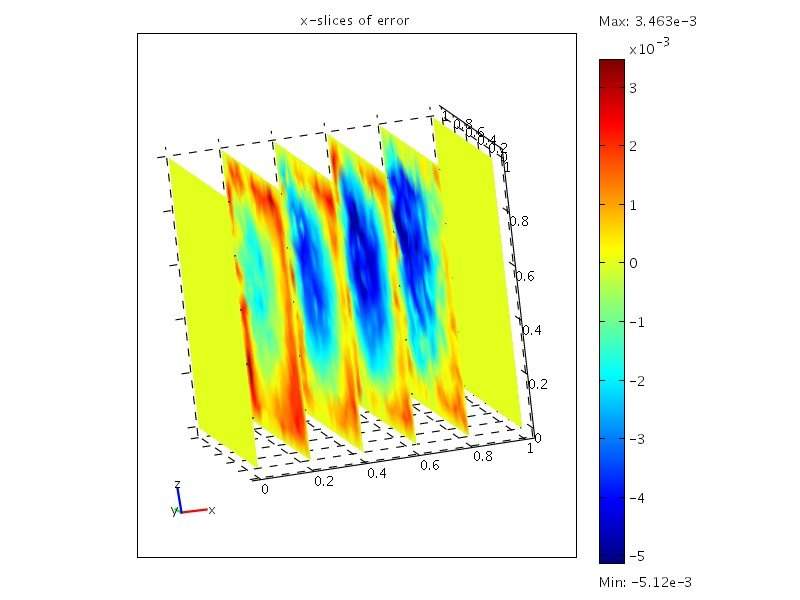}
}
\caption{{\small x-slices of the computed (moment) solution of Test 10 by
a generalized linear Hermann-Miyoshi mixed method. $\vepsi=10^{-3}$.}}\label{fig10}
\end{figure}

\begin{figure}[htb]
\centerline{
\includegraphics[angle=0,width=7cm,height=5cm]{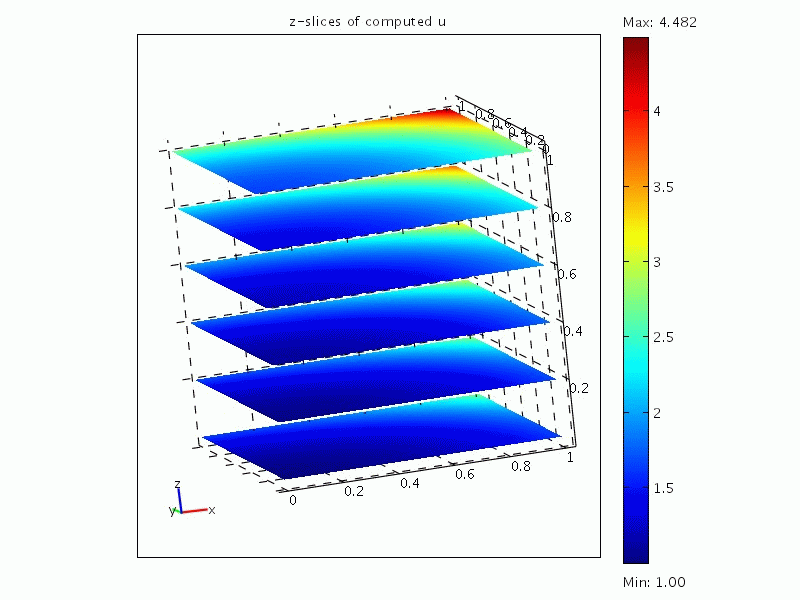}
\includegraphics[angle=0,width=7cm,height=5cm]{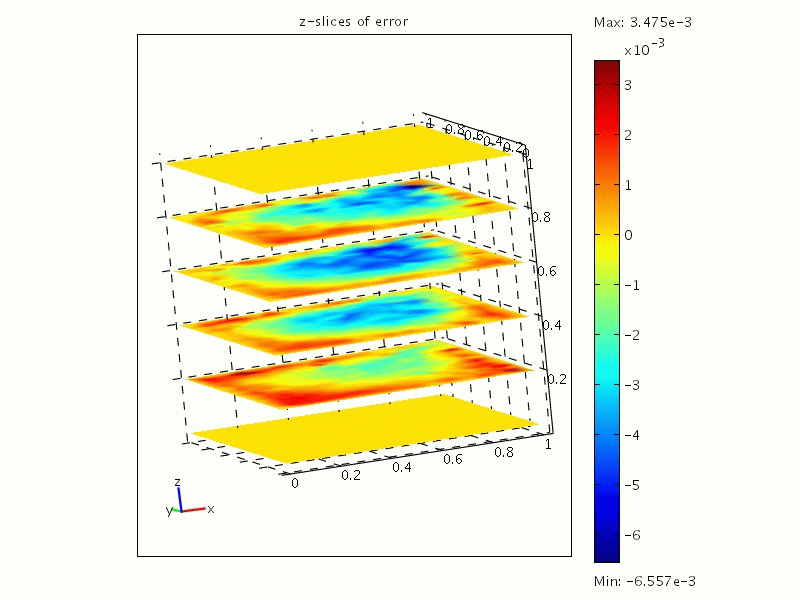}
}
\caption{{\small z-slices of the computed (moment) solution of Test 10 by
a generalized linear Hermann-Miyoshi mixed method. $\vepsi=10^{-3}$.}}\label{fig11}
\end{figure}

Figure \ref{fig10} displays color plots of five x-slices of the 
computed (moment) solution $u^\vepsi_h$ (left graph) and 
its corresponding error (right graph). 
Figure \ref{fig11} displays color plots of five z-slices 
of the computed (moment) solution $u^\vepsi_h$ (left graph) 
its corresponding error (right graph). As expected, 
the vanishing moment approximation correctly captures the 
convex viscosity solution $u^0$. Hence the
moment solution coincides with the viscosity solution.

\smallskip
{\bf Test 11:}  Our last numerical test solves the $3$-dimensional 
generalization of the test problem in Test 1.  That is, we assume 
$u$ satisfies the Monge-Amp\'ere problem \eqref{ma}--\eqref{mab} in
$\Ome=(0,1)^3$ with the data
\[
f(x,y,z)\equiv 1, \qquad g(x,y,z)\equiv 0.
\]
We remark that the above problem has a unique convex viscosity 
solution but does not have a classical solution
(cf. \cite{Gutierrez01,Dean_Glowinski06b}). There is no
explicit solution formula for the boundary value problem.

\begin{figure}[htb]
\centerline{
\includegraphics[angle=0,width=7cm,height=5cm]{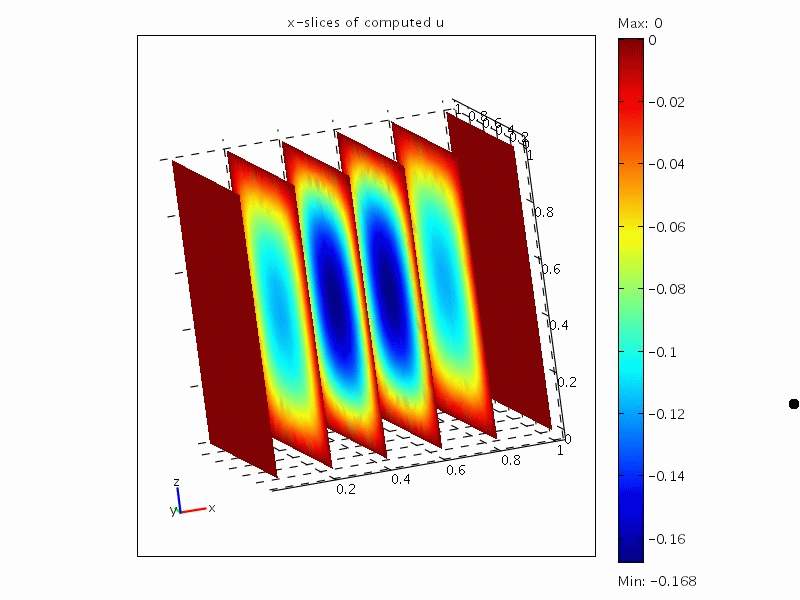}
\includegraphics[angle=0,width=7cm,height=5cm]{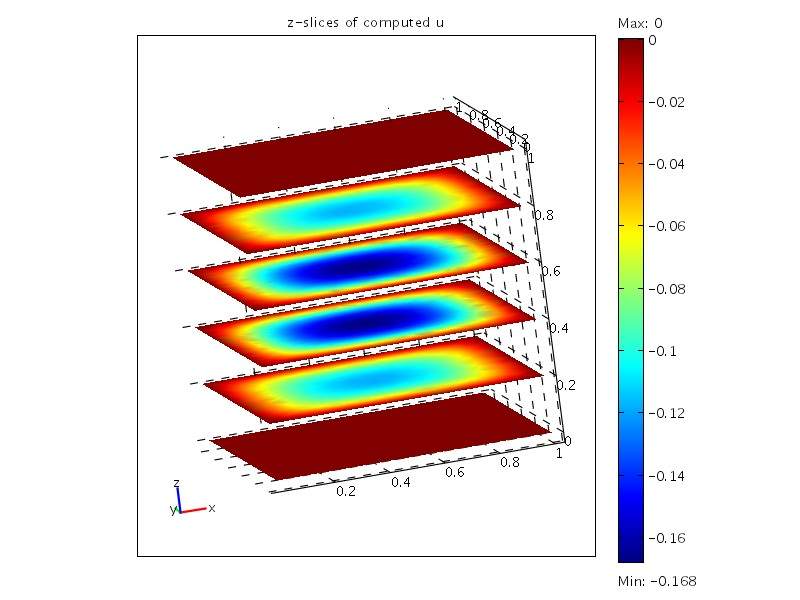}
}
\caption{{\small x-slices (left) and z-slices (right) of the 
computed (moment) solution of Test 11 by a generalized 
linear Hermann-Miyoshi mixed method. $\vepsi=10^{-3}$.}}\label{fig12}
\end{figure}

Figure \ref{fig12} displays color plots of 
x-slices (left graph) and z-slices (right graph) of the 
computed (moment) solution $u^\vepsi_h$ using
a generalized linear Hermann-Miyoshi mixed method. Once again,
the vanishing moment approximation correctly captures 
the convex viscosity solution $u^0$. Hence the
moment solution coincides with the viscosity solution.

\section{Conclusions}

In this paper we introduce a new notion of weak solutions, called 
{\em moment solutions}, through a constructive limiting 
process, called {\em the vanishing moment method}, for 
second order fully nonlinear PDEs. The notion of moment solutions
and the vanishing moment method are exactly in the same 
spirit as the original notion of viscosity solutions and 
the vanishing viscosity method proposed by M. Crandall and P. L. Lions
in \cite{Crandall_Lions83} for the Hamilton-Jacobi equations,
which is based on the idea of approximating a fully nonlinear 
PDE by a higher order quasilinear PDE. We first  present a general
framework of the vanishing moment method and the notion of 
moment solutions in \S \ref{sec-2}. We then apply the
general framework to several classes of PDEs including
the Monge-Amp\'ere type equations, Pucci's extremal equations, 
the infinite Laplace
equation, and second order fully nonlinear parabolic PDEs.
We then propose two classes of numerical methods to
discretize the fourth order ``regularized/perturbed" 
vanishing moment approximation equations. Finally, we 
present a number of numeral experiments using 
the vanishing moment methodology together with the proposed
numerical methods to demonstrate
convergence and effectiveness of the vanishing moment method,
as well as the relationship between the notion of
moment solutions and the notion of viscosity 
solution for second order fully nonlinear PDEs.

This paper provides a practical and systematic
methodology/approach, which can be backed by rigorous PDE 
and numerical theories, for approximating
second order fully nonlinear PDEs. As a by-product, 
the moment solution theory will provide some insights to our
understanding of the viscosity solution theory, and might 
provide a logical and natural generalization/extension
for the notion of viscosity solution,
especially, in the cases where there are no theories or
the existing viscosity solution theory fails (such as Monge-Amp\`ere 
equations of sub elliptic and hyperbolic types \cite{Chang_Gursky_Yang02},
and systems of second order fully nonlinear PDEs.)

{\bf Acknowledgment}: The first author would like to thank his former
colleague Professor Bo Guan of Ohio State University for many 
stimulating discussions and helpful suggestions.


\end{document}